\documentclass[final]{siamltex}

\usepackage{lineno}
\usepackage{algorithmic} 
\usepackage{algorithm}
\usepackage{float}

\usepackage{lipsum}

\usepackage{subfigure}

\usepackage{amssymb}
\usepackage{wrapfig}
\usepackage{graphicx}

\usepackage{amsmath}
\usepackage{lineno}

\renewcommand{\vec}[1]{\ensuremath\boldsymbol{#1}}
\newcommand{\norm}[1]{\left|\left|#1\right|\right|}
\newcommand{\gphi}{\nabla\phi}
\newcommand{\ngphi}{\left|\left|\nabla\phi\right|\right|}
\title{The Semi Implicit Gradient Augmented Level Set Method}

\author{Ebrahim M. Kolahdouz\footnotemark[1] \and David Salac\footnotemark[1]\ \footnotemark[2]}

\begin{document}

	\maketitle

	\renewcommand{\thefootnote}{\fnsymbol{footnote}}
	\footnotetext[1]{University at Buffalo, Department of Mechanical Engineering,	Buffalo, NY, 14260}
	\footnotetext[2]{Corresponding Author {\tt davidsal@buffalo.edu}}
	\renewcommand{\thefootnote}{\arabic{footnote}}

	\begin{abstract}
		Here a semi-implicit formulation of the gradient augmented level set method is presented. 
		By tracking both the level set and it's gradient accurate subgrid information is provided, leading to highly accurate descriptions of a moving interface. 
		The result is a hybrid Lagrangian-Eulerian method that may be easily applied in two or three dimensions.
		The new approach allows for the investigation of interfaces evolving by mean curvature and by the intrinsic Laplacian of the curvature.
		In this work the algorithm, convergence and accuracy results are presented. Several numerical experiments in both two and three dimensions demonstrate the stability of the scheme.		

	\end{abstract}


%%%%%%%%%%%%%%%%%%%%%%%%%%%%%%%%%%%%%%%%%%%%%%%%%%%%%%%%%%%%%%%%%%%%%%%%%%%%%%%%%%%%%%%%%%%%%%%%%%%%%%%%%%%%%%%%%%%%%%%%%%%
\section{Introduction and Overview}
\label{sec:1.0}
\begin{figure}[ht!]
	\centering
	\label{fig:sampleLS}
	\includegraphics[width=0.6\textwidth]{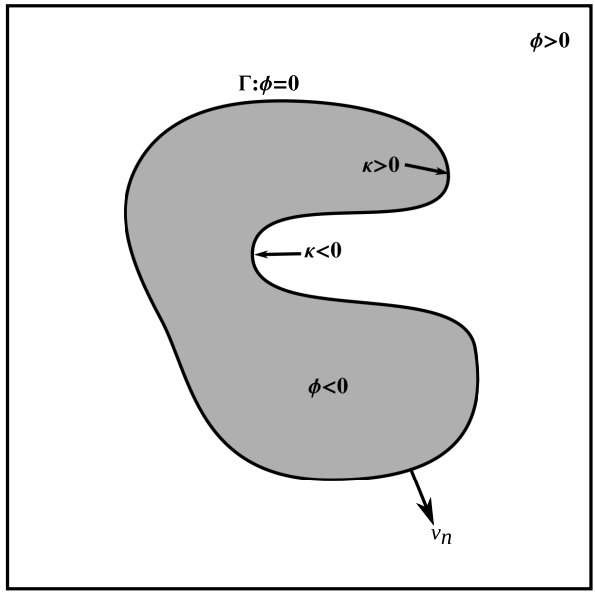}
	\caption{An example level set domain with the interface $\Gamma$, interior $\phi<0$, exterior $\phi>0$ shown. The sign of the curvature, $\kappa$,
	and the interface normal velocity $v_n$, is also shown. }	
\end{figure}

The level set method is a powerful technique to model the motion of interfaces in many disciplines. The range of applications using the level set method has 
grown substantially, including investigations of electromigration \cite{Li1999}, topology optimization \cite{Luo2008}, image processing \cite{Malladi1995},
and evolving fluid interfaces \cite{Sethian2003}. 
First developed by Osher and Sethian \cite{Osher1988} the level set method is based upon representing an interface as the zero level set of a higher dimensional function. 
Denoting the interface by $\Gamma$ and the level set function by $\phi$, we write: 
\begin{equation}
	\Gamma(t)=\left\{\vec{x}:\phi(\vec{x},t)=0\right\}.
\end{equation}
The level set value is usually chosen to be negative in the domain enclosed by interface and positive outside of it, see Fig. \ref{fig:sampleLS}. 

An implicit representation of the interface has the advantage of treating all topology changes naturally without any complex remeshing. 
Various geometric quantities such as the normal vector and curvature can be easily computed as
\begin{align}
	\vec{n}=&\frac{\gphi}{\ngphi}, \\
	\kappa=&\nabla\cdot\frac{\gphi}{\ngphi}.
\end{align}

The motion of the interface due to an underlying flow field $\vec{v}$ can be modelled by a standard advection equation,
\begin{equation}
	\label{eq:level-set-advection}
	\frac{\partial \phi}{\partial t}+\vec{v}\cdot\gphi=0.
\end{equation} 
If the normal velocity on the interface is known the equivalent advection equation is
\begin{equation}
	\label{eq:level-set-advection-normal}
	\frac{\partial \phi}{\partial t}+v_n\norm{\gphi}=0.
\end{equation} 

For accuracy reasons it is preferential to keep $\phi$ as a signed distance function such that $\ngphi=1$.
In many situations Eq. \ref{eq:level-set-advection} or Eq. \ref{eq:level-set-advection-normal} do not enforce this condition. A reinitialization process is typically used to retain the signed-distance property of $\phi$ by solving a secondary PDE in the form of \cite{Sussman1994}
\begin{equation}
\frac{\partial\phi}{\partial\tau} =\textrm{sign}(\phi_{0})\cdot(\ngphi-1),
\end{equation}
where $\tau$ is a pseudo-time taken proportional to the grid spacing.

In many systems the motion of an interface is influenced by the curvature or the variation of the curvature (surface diffusion) along an interface 
\cite{Li1999,Salac2011,Adalsteinsson,MERRIMAN1994,Mullins1957,Veerapaneni2009}. The dependence of the flow on the curvature causes numerical stability issues
for many methods. Explicit schemes require small time steps, on the order of $h^2$ for curvature based flows and $h^4$ for 
surface diffusion, where $h$ is the grid spacing \cite{chopp1999,Khenner2001} . These stringent restrictions lead to long computation times and can limit the size of the system that may be 
investigated. A fully implicit method would alleviate the time step restriction but would add the complexity of solving a possibly non-linear system every time step.
Another approach is to extract the linear portion of the advection equation and treat it implicitly, with the non-linear part treated explicitly
\cite{Salac2008,Smereka2003}. This has the advantage of numerically stabilizing the system at the cost of a slight deterioration of the accuracy of the solution.

The gradient augmented level set method is an extension of the standard level set method which advances the level set and the gradient of the level set field
\cite{Nave2010}.
This is accomplished by advecting the level set and its gradient as independent quantities in a coupled manner, ensuring that the relationship between the 
quantities remains.
This additional information allows for the determination of data on a sub-grid level, allowing for the use of coarser meshes. The original gradient augmented 
level set method was demonstrated using analytic flow fields and demonstrated excellent accuracy with minimal additional computational cost. The use of the augmented gradient
level set method in modelling vesicle motion \cite{Salac2011} has shown that while the time step restriction is not as stringent as regular explicit schemes it is not
possible to take time steps on the order of the grid spacing, as is possible with semi-implicit methods \cite{Salac2008,Smereka2003}.

In this article a semi-implicit gradient augmented level set (SIGALS) method is presented. This extension allows for the investigation of both mean curvature and surface diffusion
based motions, with restrictions based on accuracy and not stability. The tracking of the level set gradient function provides more accuracy and alleviates some of
the non-local behavior of the original semi-implicit level set method \cite{Smereka2003}. In Sec. \ref{sec:2.0} the standard gradient augmented method and the original semi-implicit method
are briefly outlined. A description of the semi-implicit gradient augmented level set method follows. Convergence analysis for a simple case and sample results are presented
in Sec. \ref{sec:3.0}. Section \ref{sec:4.0} provides a short conclusion.

%%%%%%%%%%%%%%%%%%%%%%%%%%%%%%%%%%%%%%%%%%%%%%%%%%%%%%%%%%%%%%%%%%%%%%%%%%%%%%%%%%%%%%%%%%%%%%%%%%%%%%%%%%%%%%%%%%%%%%%%%%%
\section{Development of the Method}
\label{sec:2.0}

This section will outline the standard semi-implicit level set scheme, the original augmented level set method, and the proposed new method. This section concludes 
by providing the complete algorithm of the proposed method.

%%%%%%%%%%%%%%%%%%%%%%%%%%%%%%%%%%%%%%%%%%%%%%%%%%%%%%%%%%%%%%%%%%%%%%%%%%%%%%%%%%%%%%%%%%%%%%%%%%%%%%%%%%%%%%%%%%%%%%%%%%%
\subsection{The Semi Implicit Level Set Formulation}
\label{sec:SILS}

The level set advection equation, Eq. (\ref{eq:level-set-advection}), can be written as a Hamilton-Jacobi equation,
\begin{equation}
	\frac{\partial \phi}{\partial t}+H\left(\vec{x},t,\phi,\nabla \phi\right)=0.
\end{equation}
A Hamilton-Jacobi equation can be made semi-implicit one of two ways. The first is to extract the linear portion of the Hamiltonian, $H$, and treat it implicitly,
\begin{equation}
	\label{eq:explicit-implicit}
	\frac{\partial \phi}{\partial t}+L\left(\vec{x},t^{n+1},\phi^{n+1},\nabla \phi^{n+1}\right)+N\left(\vec{x},t^{n},\phi^{n},\nabla \phi^{n}\right)=0,
\end{equation}
where $L$ is the linear portion and $N$ is the nonlinear portion such that $H=L+N$. If it is not possible to explicitly extract the linear portion one can always
determine an approximation to the linear portion and solve the following semi-implicit equation,
\begin{equation}
	\label{eq:implicit-implicit}
	\frac{\partial \phi}{\partial t}+H\left(\vec{x},t^{n},\phi^{n},\nabla \phi^{n}\right)+\widetilde{L}\left(\vec{x},t^{n},\phi^{n},\nabla \phi^{n}\right)
	-\widetilde{L}\left(\vec{x},t^{n+1},\phi^{n+1},\nabla \phi^{n+1}\right)=0,
\end{equation}
where $\widetilde{L}$ is the approximate linear portion. For mean curvature based flows it is typical to take $\widetilde{L}=\beta\nabla^2 \phi$ while for surface diffusion
based flows $\widetilde{L}=\beta\nabla^4 \phi$, where $\beta$ is a constant and typically taken to be $\beta=1/2$ \cite{Salac2008,Smereka2003}. 
This additional smoothing allows the semi-implicit level set method to utilize much larger time steps than is possible with an explicit scheme \cite{Smereka2003}.
For example, the motion by surface diffusion of a seven-lobed star requires a CFL condition of $\Delta t/h^4\approx 0.25$ using an explicit method \cite{chopp1999}
compared to a CFL condition of $\Delta t/h^4\approx 2\times 10^4$ for the semi-implicit scheme \cite{Smereka2003}.

%%%%%%%%%%%%%%%%%%%%%%%%%%%%%%%%%%%%%%%%%%%%%%%%%%%%%%%%%%%%%%%%%%%%%%%%%%%%%%%%%%%%%%%%%%%%%%%%%%%%%%%%%%%%%%%%%%%%%%%%%%%
\subsection{Standard Gradient Augmented Level Set Method}
\label{sec:GALS}

\begin{figure}[ht!]
	\centering	
	\includegraphics[width=0.4\textwidth]{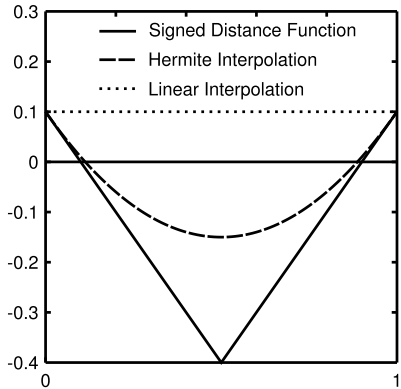}
	\caption{A one-dimensional grid with data provided at grid points $x=0$ and $x=1$. The linear interpolant misses the interface existing between the two grid points. The
	Hermite interpolant determined using the additional gradient information finds two interfaces close to the true interface locations.}
	\label{fig:sample_GALS}
\end{figure}

The gradient augmented level set method is an extension of the standard level set method which advects both the level set, $\phi$, and the level set gradient,
$\vec{\psi}=\gphi=\left(\phi_x,\phi_y,\phi_z\right)=\left(\psi^x,\psi^y,\psi^z\right)$:
\begin{align}
	\label{eq:GALS_phi}	\frac{\partial\phi}{\partial t}+\vec{v}\cdot\gphi=&0,\\
	\label{eq:GALS_psi} \frac{\partial\vec{\psi}}{\partial t}+\vec{v}\cdot\nabla \vec{\psi}+\nabla\vec{v}\cdot\vec{\psi}=&0.
\end{align}
The inclusion of gradient information allows for the determination of sub-grid information. Take as an example two grid points, $x_0=0$ and $x_1=1$, on a one
dimensional grid with $\phi\left(x_0\right)=\phi\left(x_1\right)=0.1$ with $\phi_x\left(x_0\right)=1$ and $\phi_x\left(x_1\right)=-1$. The exact level set, the 
signed distance function, gives an interface ($\phi=0$) at points $x=0.1 \textrm{ and } 0.9$. Using only the level set values linear interpolation does not return any 
interface points in this domain. The additional gradient information allows for the determination of a Hermite interpolant giving interface locations of
$x\approx 0.112$ and $x\approx 0.887$. See Fig. \ref{fig:sample_GALS} for a graphical representation.
	
To ensure that the level set and gradient field remain coupled throughout time Eqs. (\ref{eq:GALS_phi}) and (\ref{eq:GALS_psi}) are advanced in a coherent and fully coupled
manner \cite{Nave2010}. Lagrangian techniques are used to trace characteristics back in time to determine departure locations. 
Using the available information a Hermite interpolating polynomial is calculated and utilized to determine the departure values for the level set and gradient functions. To
first order this results in the following scheme,
\begin{align}
	\vec{x}_d&=\vec{x}-\Delta t\;\vec{v}\left(\vec{x},t\right), \\
	\nabla{\vec{x}}_d&=\vec{I}-\Delta t\;\nabla \vec{v}\left(\vec{x},t\right), \\
	\phi\left(\vec{x},t\right)&=P\left(\vec{x}_d,t-\Delta t\right), \\
	\vec{\psi}\left(\vec{x},t\right)&=\nabla{\vec{x}}_d\cdot \vec{G}\left(\vec{x}_d,t-\Delta t\right),
\end{align}
where $P(\vec{x},t)$ is the Hermite interpolating polynomial of $\phi$ at a time $t$, $\vec{G}=\nabla P$ is the gradient of the Hermite interpolant
defined for the level set, while $\vec{x}_d$ and $\nabla{\vec{x}}_d$ are the departure location and ``departure" gradient, respectively. Implementing a third-order version of the method Nave \textit{et. al.} demonstrated that for analytic flow fields the gradient augmented provides
results comparable to higher-order WENO schemes at a lower computational cost \cite{Nave2010}. Results for flows depending on derivatives of the level set function where not
presented.

%%%%%%%%%%%%%%%%%%%%%%%%%%%%%%%%%%%%%%%%%%%%%%%%%%%%%%%%%%%%%%%%%%%%%%%%%%%%%%%%%%%%%%%%%%%%%%%%%%%%%%%%%%%%%%%%%%%%%%%%%%%
\subsection{The Semi-Implicit Gradient Augmented Level Set Method}
\label{sec:2.3}

The semi-implicit gradient augmented level set method is a combination of the methods briefly explained in Secs. \ref{sec:SILS} and \ref{sec:GALS}. The goal is to 
combine the additional accuracy afforded by explicitly tracking gradient information with the stability properties of a semi-implicit scheme.
Equations (\ref{eq:GALS_phi}) and (\ref{eq:GALS_psi}) are replaced with 
\begin{align}
	\label{eq:SIGALS_phi} \frac{D\phi}{Dt}+\beta L\phi-\beta L\phi&=0, \\
	\label{eq:SIGALS_psi} \frac{D\vec{\psi}}{Dt}+\nabla\vec{v}\cdot\vec{\psi}+\beta L\vec{\psi}-\beta L\vec{\psi}&=0,
\end{align}
where $D/Dt$ is the material (Lagrangian) derivative, $L$ is a linear operator, and $\beta$ is a constant.
To first-order in time this can be written as
\begin{align}
	\label{eq:departure} \vec{x}_d&=\vec{x}-\Delta t\;\vec{v}\left(\vec{x},t\right), \\
	\label{eq:departure_grad} \nabla{\vec{x}}_d&=\vec{I}-\Delta t\;\nabla \vec{v}\left(\vec{x},t\right), \\
	\label{eq:phih} \widetilde{\phi}&=P\left(\vec{x}_d,t-\Delta t\right), \\
	\label{eq:psih} \widetilde{\vec{\psi}}&=\nabla{\vec{x}}_d\cdot \vec{G}\left(\vec{x}_d,t-\Delta t\right), \\
	\label{eq:Iphi} \frac{\phi_{n+1}-\widetilde{\phi}}{\Delta t}&=\beta L \phi_{n+1} - \beta L \phi_n, \\
	\label{eq:Ipsi} \frac{\vec{\psi}_{n+1}-\widetilde{\vec{\psi}}}{\Delta t}&=\beta L \vec{\psi}_{n+1} - \beta L \vec{\psi}_n, 
\end{align}
The linear operator $L$ is based on the underlying flow field. For mean curvature flow $L=\nabla^2$ while for surface diffusion 
$L=\nabla^4$. If $\beta=0$ 
the method results in the standard gradient augmented method. If instead of Eqs. \ref{eq:departure} to \ref{eq:psih} values are set as 
$\widetilde{\phi}=\phi_n$ 
and $\widetilde{\vec{\psi}}=\vec{\psi}_n$ the method results in the standard semi-implicit level set scheme.

The location $\vec{x}_d$ is obtained by tracing characteristics backwards in time. In general
this location will not lie on a grid point and thus the use of an interpolant, $P$ and $\vec{G}=\nabla P$, is required. 
In two dimensions let $\vec{x}_d$ lie within a grid cell $\Omega_{i,j}$ enclosing the region given by four grid points:  $\vec{x}_{i,j}$, $\vec{x}_{i+1,j}$, 
$\vec{x}_{i,j+1}$, and $\vec{x}_{i+1,j+1}$. Using the value of level set, $\phi$, and gradient field, $\vec{\psi}=(\psi^x,\psi^y)$,
at time $t-\Delta t$ it is possible to define the 
Hermite interpolant over the grid cell by requiring that 
$P(\vec{x}_{m,n})=\phi_{m,n}$, $\partial_x P(\vec{x}_{m,n})=\psi^x_{m,n}$, $\partial_y P(\vec{x}_{m,n})=\psi^y_{m,n}$, and 
$\partial_{xy} P(\vec{x}_{m,n})=(\partial_x \psi^y_{m,n}+\partial_y \psi^x_{m,n})/2$ for $m={i,i+1}$ and $n={j,j+1}$. The gradient interpolant is then
defined as $\vec{G}=\nabla P$.

Two velocity fields are considered here. The first velocity field is mean curvature flow given by $v_n=-\kappa$ or $\vec{v}=-\kappa\vec{n}$, where $\kappa$ is the mean curvature and
$\vec{n}$ is the unit normal vector. The unit normal vector is simply given as the normalized gradient vector:
\begin{equation}
	\vec{n}=\frac{\vec{\psi}}{\|\vec{\psi}\|}.
\end{equation}
The curvature is computed as
\begin{equation}
	\label{eq:2DK}
	\kappa=\frac{\phi_{xx}\phi_y^2+\phi_{yy}\phi_x^2-2\phi_{xy}\phi_x\phi_y}{\left(\phi_x^2+\phi_y^2+\epsilon\right)^{3/2}}
\end{equation}
for two-dimensional flows and 
\begin{align}
	\label{eq:3DK}
	\kappa=&(\phi_{xx}\left(\phi_y^2+\phi_z^2\right)+\phi_{yy}\left(\phi_x^2+\phi_z^2\right)+\phi_{zz}\left(\phi_x^2+\phi_y^2\right) \nonumber \\
	&-2\phi_{xy}\phi_x\phi_y-2\phi_{xz}\phi_x\phi_x-2\phi_{yz}\phi_y\phi_z)/\left(\phi_x^2+\phi_y^2+\phi_z^2+\epsilon\right)^{3/2}
\end{align}
for three-dimensional flows.
In the above the first-order derivatives, $\phi_x$, $\phi_y$, and $\phi_z$ are replaced by their corresponding values in the $\vec{\psi}$ vector. The second-order 
derivatives are obtained by first derivatives of the $\vec{\psi}$ vector. 
Any cross derivatives are averages of the two possible first derivatives of the gradient vector, \textit{i.e.} $\phi_{xy}=(\partial_x\psi^y+\partial_y\psi^x)/2$.
The value $\epsilon=10^{-8}$ is added to ensure that no division by zero will occur. 

The second velocity field considered is that of surface diffusion, given by $v_n=\nabla^2_s\kappa$ or $\vec{v}=\left(\nabla^2_s\kappa\right)\vec{n}$ where $\nabla^2_s$ is 
the surface Laplacian. Assume that the curvature is known in the vicinity of the interface. In two dimensions the surface Laplacian of the curvature can then be calculated as
\begin{equation}
	\label{eq:2DS}
	S=\nabla_s^2\kappa=\frac{\kappa_{xx}\phi_y^2+\kappa_{yy}\phi_x^2-2\kappa_{xy}\phi_x\phi_y}{\phi_x^2+\phi_y^2+\epsilon}-\kappa\frac{\kappa_x\phi_x+\kappa_y\phi_y}{\sqrt{\phi_x^2+\phi_y^2+\epsilon}},
\end{equation}
while for the three-dimensional case 
\begin{align}
	\label{eq:3DS}
	S=\nabla_s^2\kappa=&(\kappa_{xx}\left(\phi_y^2+\phi_z^2\right)+\kappa_{yy}\left(\phi_x^2+\phi_z^2\right)+\kappa_{zz}\left(\phi_x^2+\phi_y^2\right) \nonumber \\
	&-2\kappa_{xy}\phi_x\phi_y-2\kappa_{xz}\phi_x\phi_z-2\kappa_{yz}\phi_y\phi_z)/(\phi_x^2+\phi_y^2+\phi_z^2+\epsilon) \nonumber \\
	&-\kappa ( \kappa_x\phi_x+\kappa_y\phi_y+\kappa_z\phi_z)/\sqrt{\phi_x^2+\phi_y^2+\epsilon}.
\end{align}

The curvature and surface diffusion velocity fields defined above are only defined on the interface. To determine a smooth velocity field elsewhere in the domain values can be extended from the interface through the use of 
an extension equation applied to a quantity $q$ \cite{Salac2008,Smereka2003,Peng1999},
\begin{equation}
	\label{eq:extension} \frac{\partial q}{\partial \tau}+\textrm{sign}(\phi)\vec{n}\cdot\nabla q=0,
\end{equation}
where $\tau$ is a fictitious time. 

The semi-implicit gradient augmented level set scheme can be summarized in the following two algorithms:

\begin{algorithm}                 
	\caption{Compute the Velocity Field}          
	\label{alg:vel}                           
	\begin{algorithmic}                    
		\REQUIRE $\phi$ and $\vec{\psi}$		
		\STATE Set $\kappa=0$ in the entire domain.
		\STATE Compute $\kappa$ using either Eq. (\ref{eq:2DK}) or (\ref{eq:3DK}) at grid points next to the interface.
		\STATE Extend $\kappa$ by at least 4 grid points using an extension algorithm \cite{Salac2008,Smereka2003,Peng1999}.
		\IF{velocity field is surface diffusion}
			\STATE Set $S=0$ in the entire domain.
			\STATE Compute $S$ using either Eq. (\ref{eq:2DS}) or (\ref{eq:3DS}) at grid points next to the interface.
			\STATE Extend $S$ by at least 4 grid points using an extension algorithm \cite{Salac2008,Smereka2003,Peng1999}.
			\STATE Set the velocity as $\vec{v}=S\vec{n}$.
		\ELSE
			\STATE Set the velocity as $\vec{v}=-\kappa\vec{n}$.
		\ENDIF
		\RETURN $\vec{v}$		
	\end{algorithmic}
\end{algorithm}

\begin{algorithm}                  % enter the algorithm environment
	\caption{Advance the Level Set and Gradient Fields}          % give the algorithm a caption
	\label{alg:SIGALS}                           % and a label for \ref{} commands later in the document
	\begin{algorithmic}                    % enter the algorithmic environment
		\REQUIRE $\phi_n$ and $\vec{\psi}_n$.
		\STATE Compute $\vec{v}$ using Algorithm \ref{alg:vel}.
		\FOR{every grid point in the domain}
			\STATE Compute the departure point $\vec{x}_d$ using Eq. \ref{eq:departure}.
			\STATE Compute the departure gradient $\nabla\vec{x}_d$ using Eq. \ref{eq:departure_grad}
			\STATE Evaluate the tentative Lagrangian solutions $\widetilde{\phi}$ and $\widetilde{\vec{\psi}}$ by Eqs. (\ref{eq:phih}) and (\ref{eq:psih}).			
		\ENDFOR
		\IF{velocity field is surface diffusion}
			\STATE Set $L=\nabla^4$
		\ELSE
			\STATE Set $L=\nabla^2$
		\ENDIF
		\STATE Solve $\left(\vec{I}-\Delta t \beta L \right)\phi_{n+1}=\widetilde{\phi}-\Delta t \beta L \phi_{n}$.
		\STATE Solve $\left(\vec{I}-\Delta t \beta L \right)\vec{\psi}_{n+1}=\widetilde{\vec{\psi}}-\Delta t \beta L \vec{\psi}_{n}$.
		\RETURN $\phi_{n+1}$ and $\vec{\psi}_{n+1}$.
	\end{algorithmic}
\end{algorithm}

\clearpage

%%%%%%%%%%%%%%%%%%%%%%%%%%%%%%%%%%%%%%%%%%%%%%%%%%%%%%%%%%%%%%%%%%%%%%%%%%%%%%%%%%%%%%%%%%%%%%%%%%%%%%%%%%%%%%%%%%%%%%%%%%%
\section{Numerical Results}
\label{sec:3.0}

To avoid any anisotropy introduced by the use of center finite difference approximations all derivatives will be 
approximated by isotropic finite differences \cite{Kumar2004}.
The linear systems in Eqs. (\ref{eq:Iphi}) and (\ref{eq:Ipsi}) are solved using a standard Bi-CGSTAB method \cite{Vorst1992,Saad2003}. 
All of the results will be 
first order in time and utilize second-order isotropic finite differences. In each case 
periodic boundary conditions are assumed and a constant of $\beta=0.5$ is used. The curvature and surface Laplacian of curvature are 
first calculated near the interface and then extended to the rest of the domain,
as shown in the Alg. \ref{alg:vel}. The time step is given by $\Delta t$ while the uniform grid spacing is given by $h$.
Every initial interface is initially described by a signed distance function.
Please note that no level set reinitialization was performed during the simulations.

\subsection{Mean Curvature Flow}
\label{sec:3.1}

Here the evolution of two- and three-dimensional interfaces under mean-curvature flow will be shown. For all cases in this section the 
velocity of the interface is given by $\vec{v}=-\kappa\vec{n}$.

\begin{figure}
	\centering
	\includegraphics[width=0.5\textwidth]{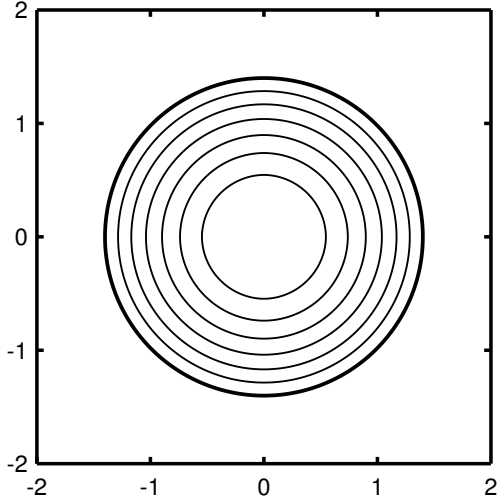}
	\caption{A circle collapsing due to mean curvature flow. The initial interface is represented by the thick line.
		The the grid spacing is $h=0.0625$ and the time step is $\Delta t=4h^2=0.015625$.
		The interface is shown in time increments of $10\Delta t$.}
	\label{fig:circleK}	
\end{figure}

First consider a circular interface. 
With mean curvature flow a circular interface beginning with an initial radius of $r_0$ will collapse uniformly.
A sample result is presented in Fig. \ref{fig:circleK} for a grid spacing of $h=0.0625$ and a time step of $\Delta t=4h^2=0.015625$.

The collapse of a circular interface under mean curvature flow is a situation with a known analytic solution and thus allows 
for the investigation of the accuracy of the SIGALS scheme.
At time $t$ the radius of the circle is given by
\begin{equation}
	r(t)=\sqrt{r_{0}^{2}-2t},
	\label{eq:exact_circle_radius}
\end{equation}
where $r_0$ is the initial radius. The interface was allowed to evolve until a time of $t=0.375$ with various grid spacings
and a time step set to $\Delta t=8h^2$. The level set and gradient field values for grid points next to the interface are compared
to a signed distance function describing a circle with the radius given in Eq. (\ref{eq:exact_circle_radius}).
The resulting errors are shown in Tables \ref{table:phi_error} and \ref{table:psi_error}.
The level set retains the second-order in space accuracy of the underlying discretization while the gradient is one-half order lower.

\begin{table}
	\centering	
	\caption{ Level set error for a circle with $v_n=-\kappa$ using $\Delta t=8h^2$.}
	\label{table:phi_error}	
	\begin{tabular}{r r r r r r}

		N & $h$ & $L_2$ & Order & $L_\infty$ & Order \\ [1ex]
		\hline\hline \\
		64 & 0.0625 		& $2.60\times 10^{-3}$ &      & $4.23\times 10^{-2}$ & \\
		128 & 0.03125 		& $6.47\times 10^{-4}$ & 2.0  & $1.37\times 10^{-2}$ & 1.63\\
		256 & 0.015625		& $1.64\times 10^{-4}$ &	1.98 & $4.8\times 10^{-3}$	& 1.51 \\
		512 & 0.0078125 	& $4.20\times 10^{-5}$ & 1.96 & $1.61\times 10^{-3}$ & 1.58\\ [1ex]
		\hline	
	\end{tabular}	

\end{table}

\begin{table}
	\centering
	\caption{ Gradient error for a circle with $v_n=-\kappa$ using $\Delta t=8h^2$.}
	\label{table:psi_error}	
	\begin{tabular}{r r r r r r}

		N & $h$ & $L_2$ & Order & $L_\infty$ & Order \\ [1ex]
		\hline\hline \\
		 64 & 0.0625		 & $6\times 10^{-4}$ &	  & $1.21\times 10^{-2}$ & \\
		 128 & 0.03125 	 & $1.68\times 10^{-4}$ & 1.84 & $3.7\times 10^{-3}$ & 1.71\\
		 256 & 0.015625	 & $5.64\times 10^{-5}$ & 1.57 & $1.7\times 10^{-3}$  & 1.12 \\
		 512 & 0.0078125   & $2.01\times 10^{-5}$ & 1.49 & $8.2\times 10^{-4}$ & 1.05\\ [1ex]
		\hline	
	\end{tabular}		
\end{table}

\begin{figure}
	\begin{center}
		\includegraphics[width=0.6\textwidth]{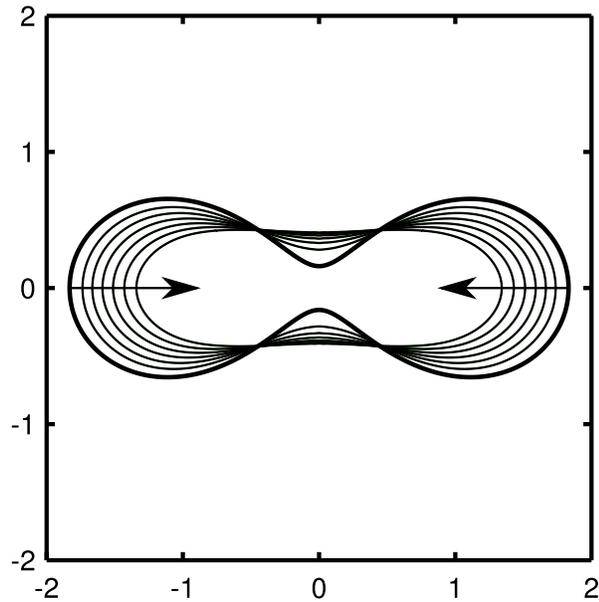}
	\end{center}
	\caption{Motion by mean curvature for a Cassini oval at different times on a $ 128\times 128$ grid. The grid spacing is $h=0.03125$ while the time step is $\Delta t=0.03$.
	The thick line represents the initial interface while the arrow indicates the direction of motion.
	Subsequent interface locations are shown in time increments of $3\Delta t$. Eventually the interface will collapse to a point.}
	\label{fig:2d_cassini_mean_curvature} 
\end{figure}

Similar behavior is observed in more complex interfaces. The collapse of a Cassini oval is shown in Fig. \ref{fig:2d_cassini_mean_curvature}.
Initially the interface has regions of both positive and negative curvature. The positive curvature regions move towards the center of the 
domain while those with negative curvature move away from the center. After some time an ellipse-like interface is obtained. From this point
forward the interface collapses to a point.

\begin{figure}
		\begin{center}
			\includegraphics[width=0.6\textwidth]{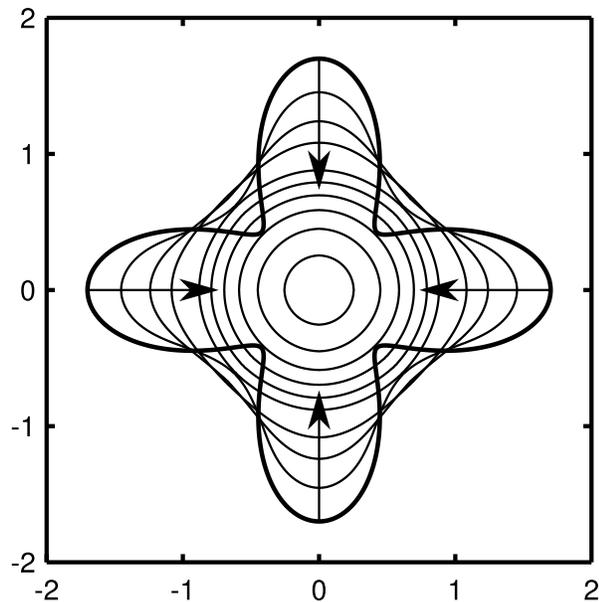}
		\end{center}
\caption{Evolution of a 4-lobe star under mean curvature flow. The grid size  $ 256\times 256 $ giving a grid spacing of $h=0.015625$.
The time step is $\Delta t=0.005$. The thick line is the initial interface with subsequent 
interface locations shown in increments of $20\Delta t$. The arrows indicate the direction of motion.}
\label{fig:2d_four_lobe_mean_curvature} 
\end{figure}

Similar results are seen for a two-dimensional four-lobe star, Fig. \ref{fig:2d_four_lobe_mean_curvature}. The interface evolves until
the curvature is strictly positive at which point the shape collapses to a point.

\begin{figure}
	\begin{center}
		\begin{tabular}[t]{cc}
			\begin{minipage}{0.35\textwidth}
				\begin{center}
					\includegraphics[width=\textwidth ]{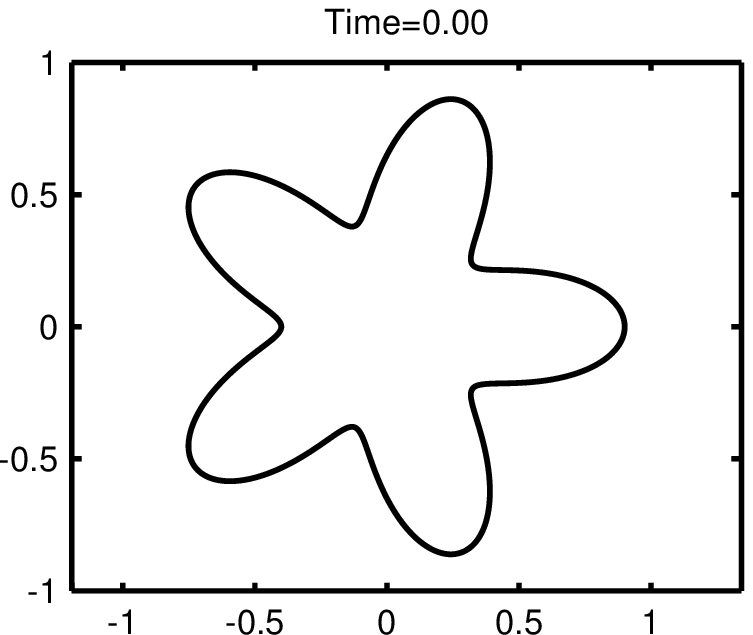} \vspace{1ex}
					
					\includegraphics[width=\textwidth ]{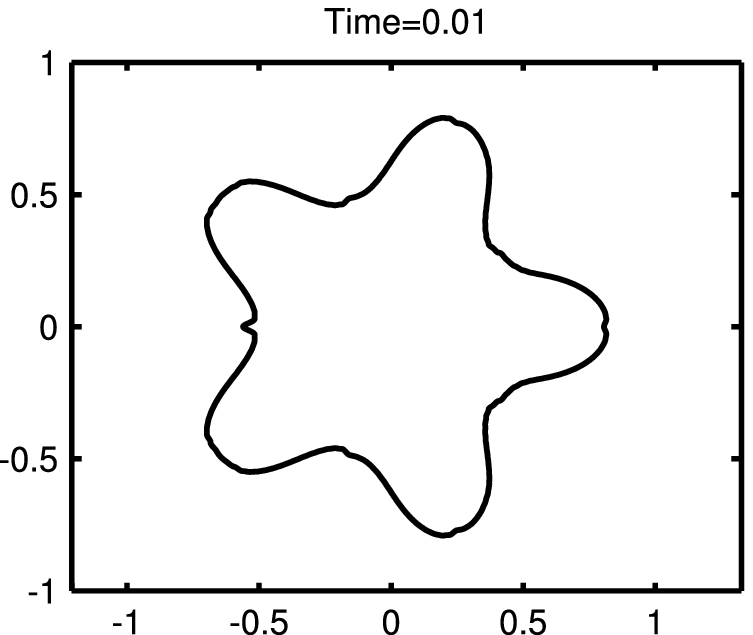}  \vspace{1ex}
					
					\includegraphics[width=\textwidth ]{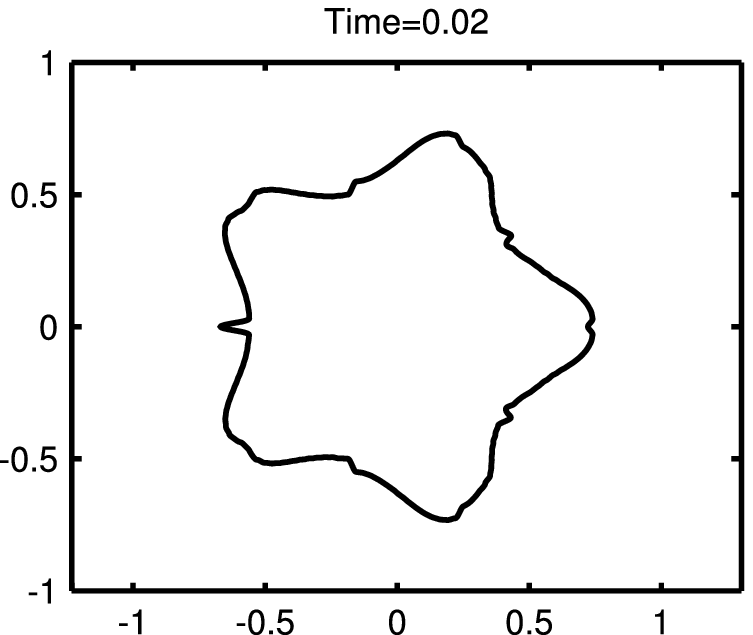}  \vspace{1ex}
					
					\includegraphics[width=\textwidth ]{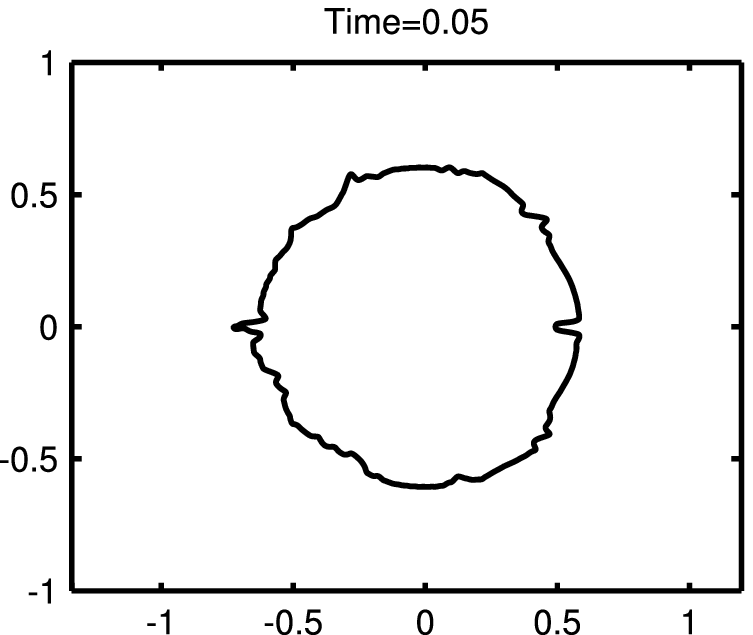}  
					
					\small Standard Gradient Augmented Level Set Method
				\end{center}
			\end{minipage} &
			\begin{minipage}{0.35\textwidth}
				\begin{center}
					\includegraphics[width=\textwidth]{initial_curvature_star}  \vspace{1ex}
					
					\includegraphics[width=\textwidth]{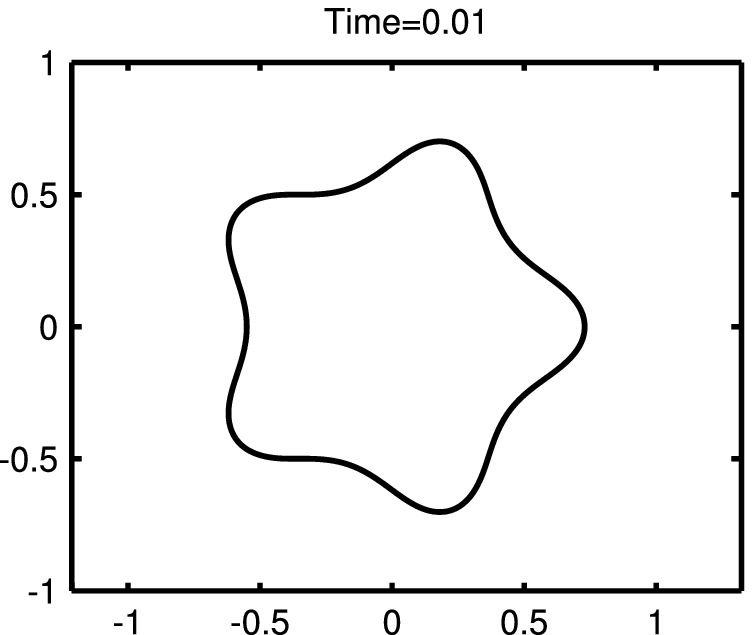}  \vspace{1ex}
					
					\includegraphics[width=\textwidth]{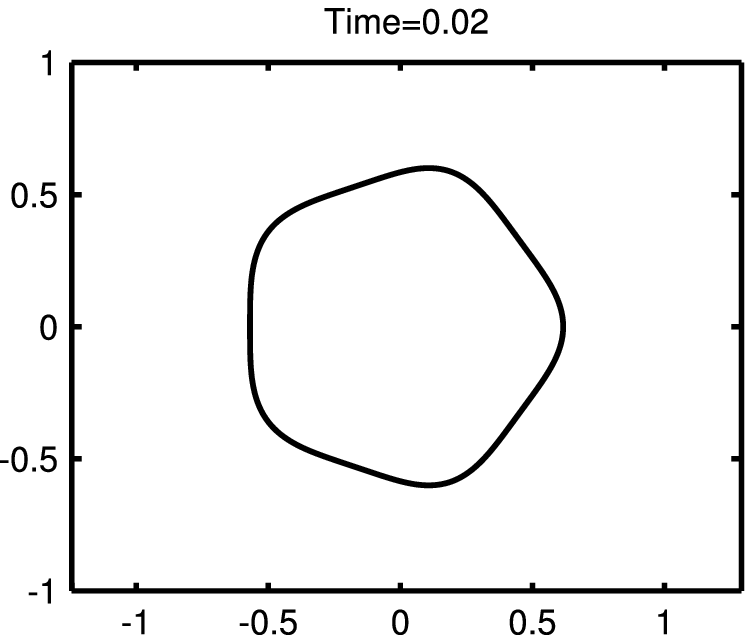}  \vspace{1ex}
					
					\includegraphics[width=\textwidth]{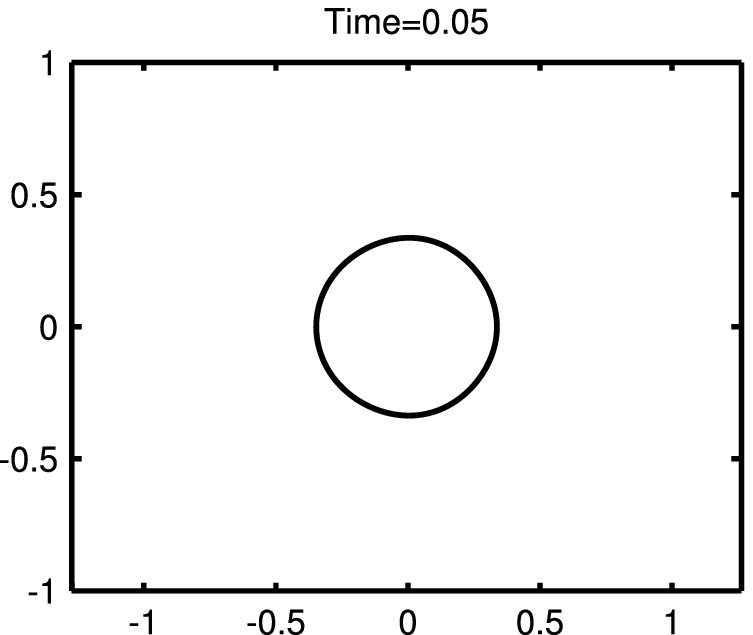} 
					
					\small Semi-Implicit Gradient Augmented Level Set Method
				\end{center}
			\end{minipage}			
		\end{tabular}								
	\end{center}
	\caption{A comparison of the standard gradient augmented level set method with the semi-implicit augmented method 
				using the same time step for both methods.
				The five-lobe star is collapsing under mean curvature flow, $v_n=-\kappa$, 
				for $\Delta t=5\times10^{-4}$ and $h=0.03125$ giving $\Delta t/h^2 \approx 0.5$. The full
				domain is the region given by $[-2,2]^2$}	
	\label{fig:2d_five_lobe_comparison}
\end{figure}

To demonstrate the added stability properties of the SIGALS scheme as compared to the original gradient augmented level set method
the evolution of a five-lobe star is shown in Fig. \ref{fig:2d_five_lobe_comparison}.
A common grid spacing of $h=0.03125$ and time step $\Delta t=5\times10^{-4}$ is used for both cases. This results in a CFL condition
of $\Delta t/h^2\approx 0.5$. Over time the standard gradient augmented method begins to demonstrate numerical instabilities. 
The SIGALS method does not demonstrate any instability and shows the expected behavior.

\begin{figure}
	\begin{center}
		\hfill
		\subfigure[$t=0$]{
			\includegraphics[width=0.475\textwidth]{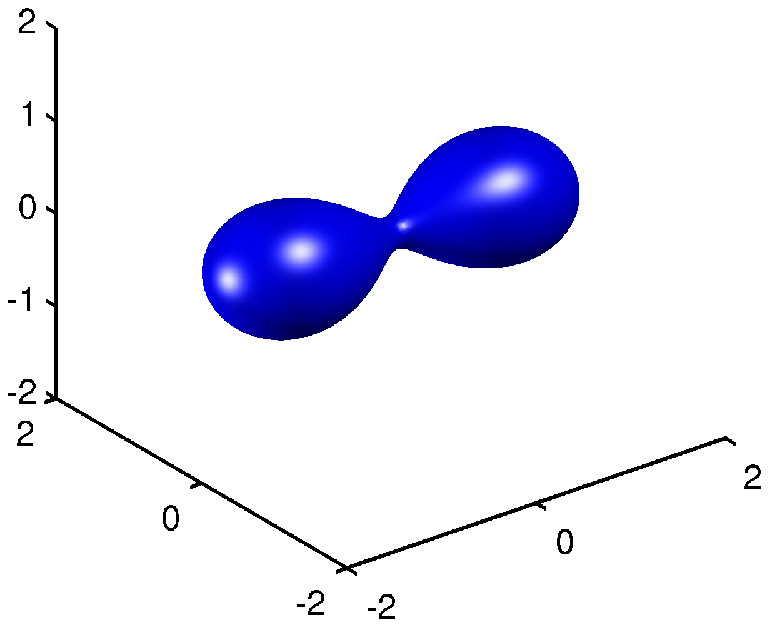}
			\label{fig:3.2.3a}
		}\hfill
		\subfigure[$t=0.025$]{
			\includegraphics[width=0.475\textwidth]{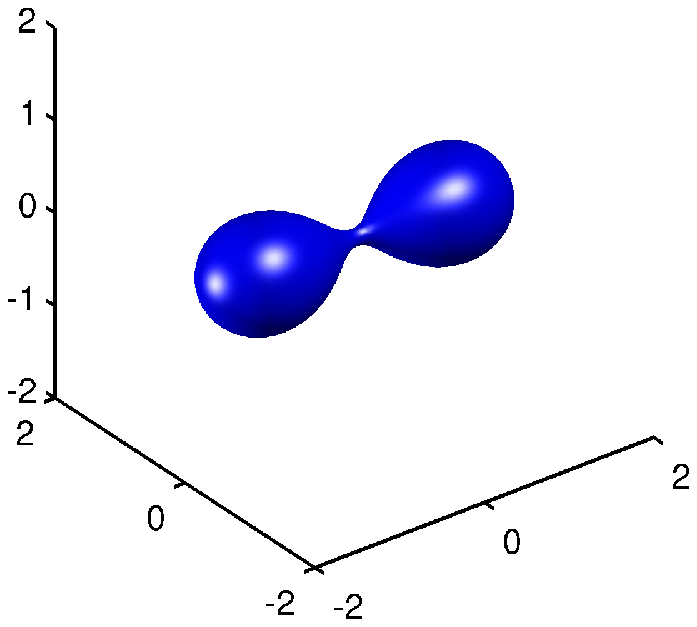}
			\label{fig:3.2.3b}
		} \hfill \\
		\hfill
		\subfigure [$t=0.04$]{
			\includegraphics[width=0.475\textwidth]{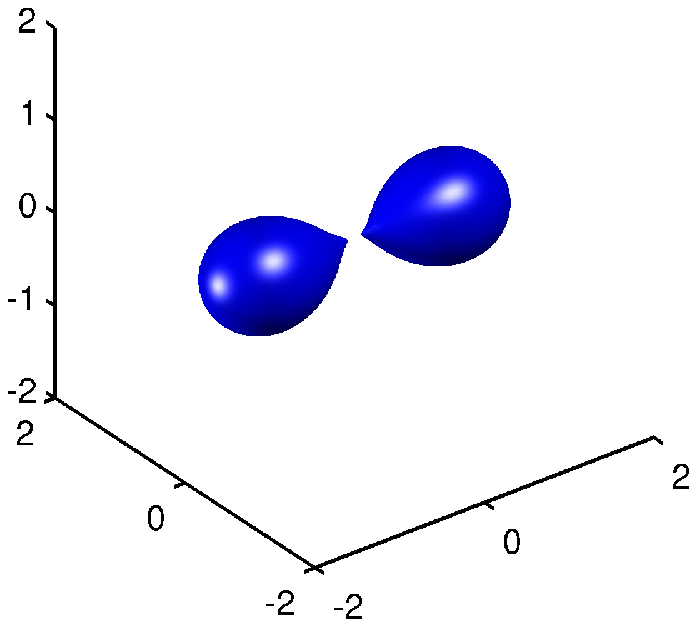}
			\label{fig:3.2.3c}
		}\hfill
		\subfigure [$t=0.075$]{
			\includegraphics[width=0.475\textwidth]{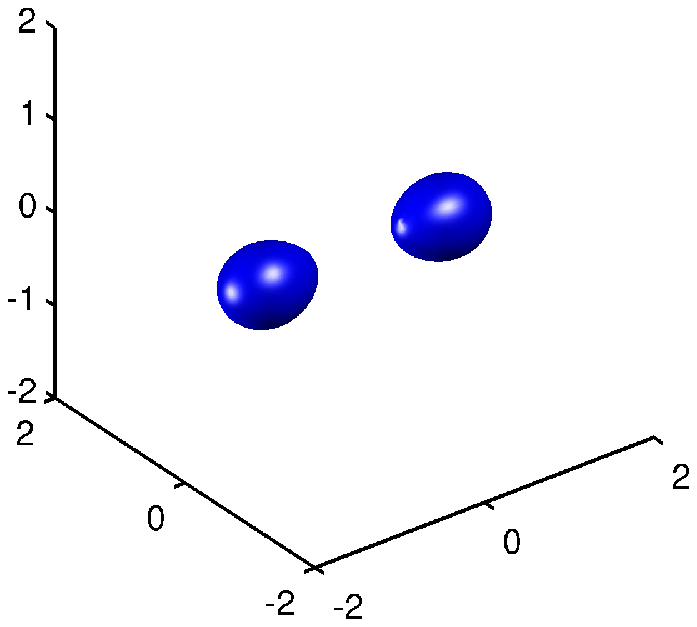}
			\label{fig:3.2.3d}
		}\hfill \\ 
		   		
		\caption{Collapse of a three dimensional Cassini oval under mean curvature flow for $h=0.0625$ and $\Delta t=0.005$.}	
		\label{fig:3d_cassini_mean_curvature}
	\end{center}				
\end{figure}

The extension of the SIGALS method to three dimensions is demonstrated in Fig. \ref{fig:3d_cassini_mean_curvature}, which shows the collapse
of a three-dimensional Cassini oval with mean-curvature flow. Due to the additional curvature the neck region does not thicken as in the 
two-dimensional Cassini oval case, instead the surface collapses with the surface splitting into two separate interfaces.
This result clearly demonstrates the level sets ability to handle
topological changes naturally.

%%%%%%%%%%%%%%%%%%%%%%%%%%%%%%%%%%%%%%%%%%%%%%%%%%%%%%%%%%%%%%%%%%%%%%%%%%%%%%%%%%%%%%%%%%%%%%%%%%%%%%%%%%%%%%%%%%%%%%%%%%%%%
\subsection{Surface Diffusion}
\label{sec:3.2}

	This section considers the motion of interfaces due to the intrinsic variation of the curvature along the interface, or simply called surface diffusion. In this case the velocity of 
	an interface is given by $\vec{v}=\left(\nabla_s^2\kappa\right)\vec{n}$, where $\nabla_s^2$ is the surface Laplacian. The final result for all surface diffusion cases will be	a constant curvature interface: a circle in two dimensions and a sphere in three dimensions.
	
\begin{figure}
		\begin{center}
			\includegraphics[width=0.6\textwidth]{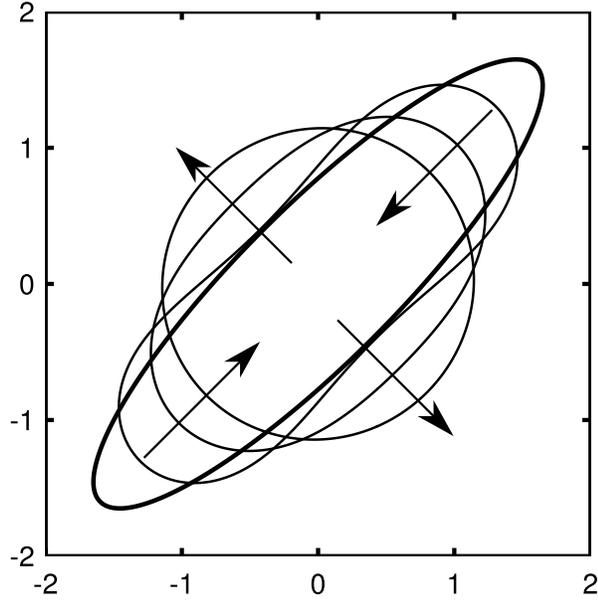}			
		\end{center}
		\caption{Motion by surface diffusion for an elliptical interface for a grid spacing of $h=0.03125$ and a time step of $\Delta t=0.001\approx 1048\Delta x^{4} $.
		The thick line is the initial interface and the arrows indicate the direction of motion. Interfaces are shown for times of 0, 0.05, 0.2, and 0.6. }
		\label{fig:2d_ellipse_diffusion} 
	\end{figure}

\begin{figure}
		\begin{center}

			\subfigure[Spatial Study]{
  				\includegraphics[width=0.475\textwidth]{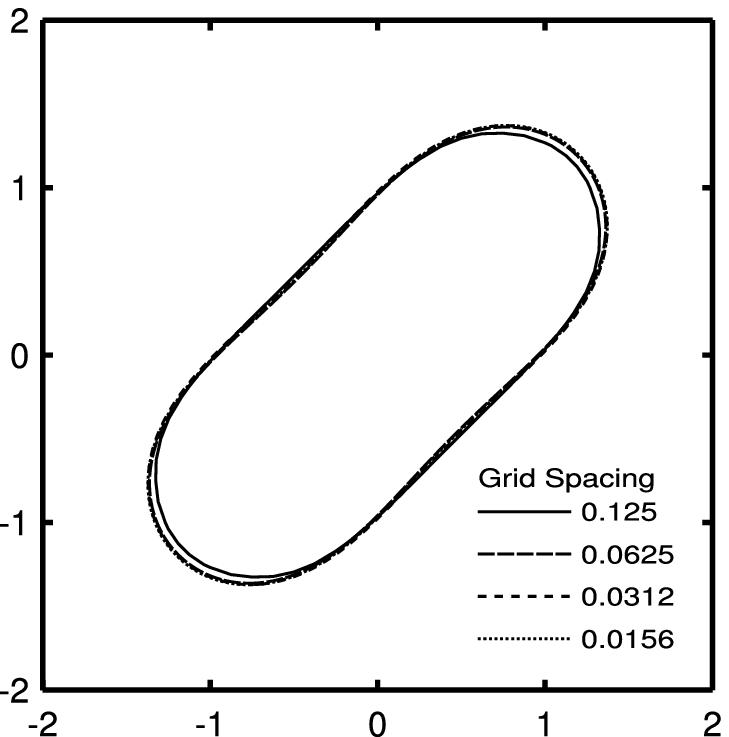}
				\label{fig:4.2.4a}
			}
			\subfigure[Time Study]{
  				\includegraphics[width=0.475\textwidth]{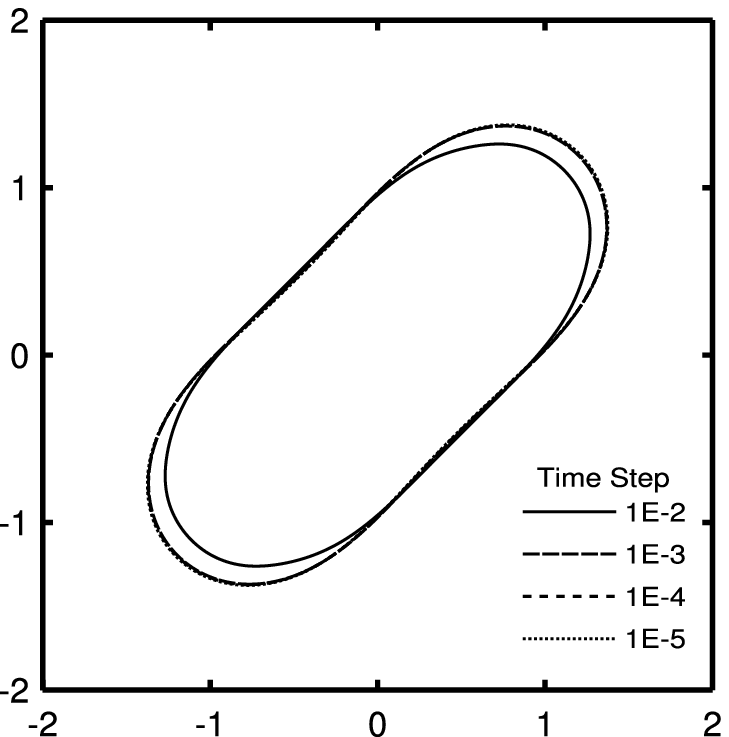}
				\label{fig:4.2.4b}
			} 
  		
			\caption{Qualitative convergence check for motion due to surface diffusion. For the spatial study the time step is fixed at $\Delta t=10^{-4}$ 
			while for the time study
			the grid spacing is fixed at $h=0.015625$. All results are shown at a time of $t=0.1$.}		
			\label{fig:2d_ellipse_convergence}
		\end{center}				
	\end{figure}

	The first shape considered is that of an inclined ellipse, Fig. \ref{fig:2d_ellipse_diffusion}.
	As expected the ellipse smoothly transitions to the circular interface. This interface is also utilized as a qualitative check on the 
	spatial and temporal convergence of the method. 
	The grid study is performed by fixing the time step at $\Delta t=10^{-4}$ and varying the grid spacing. The temporal study fixes the grid size at 
	$\Delta x=0.0312$ and varies the time step. The results for a time of 0.1 are seen in Fig. \ref{fig:2d_ellipse_convergence} and
	demonstrate the method quickly converging to a common solution.

\begin{figure}
		\begin{center}
			\subfigure[]{\includegraphics[width=0.475\textwidth]{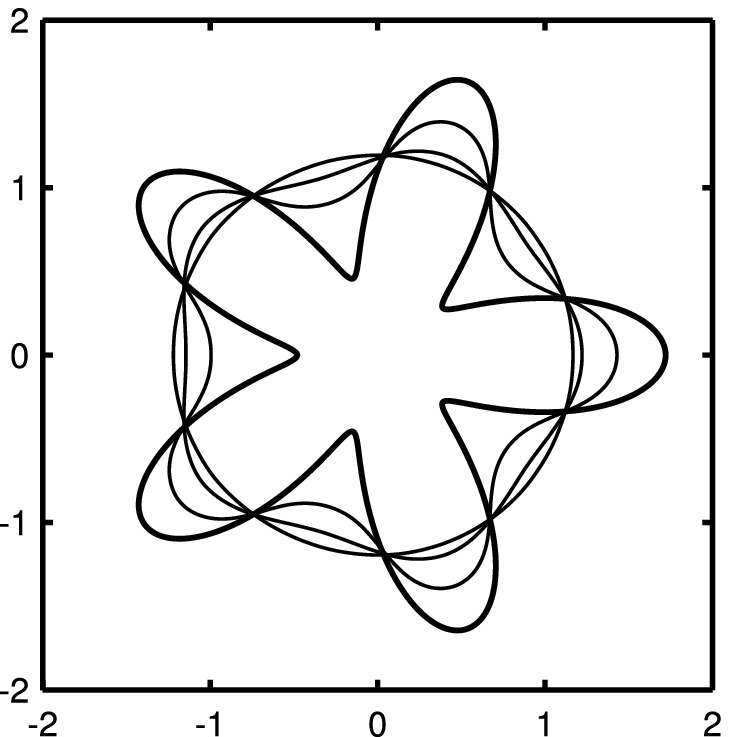}}
			\subfigure[]{\includegraphics[width=0.475\textwidth]{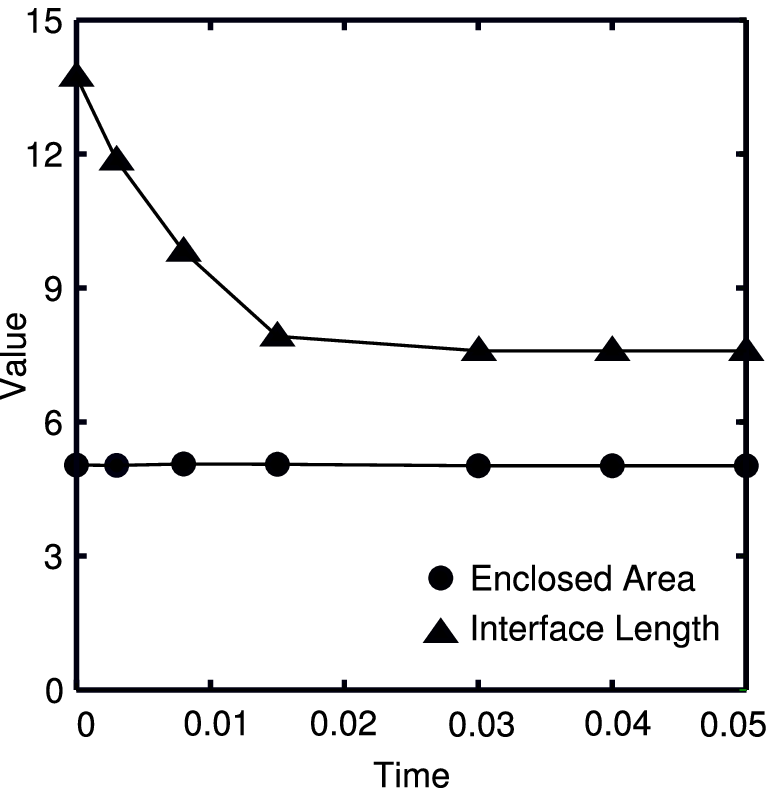}}
		\end{center}
		\caption{Motion by surface diffusion for a five-lobe star at using a grid spacing $\Delta x=0.0625 $ and a time step $\Delta t=0.001\approx 65\Delta x^{4}$. The motion
		of the interface is seen in (a) while the enclosed area and interface length are tracked in (b).}
		\label{fig:2d_star_diffusion} 
	\end{figure}

	The evolution of a five-lobed star under surface diffusion is seen in Fig. \ref{fig:2d_star_diffusion}. 
	In addition to the location of the interface over time the 
	total enclosed area and interfacial length are also presented. 
	It is clear that the enclosed area remains constant at $\pi$ and the interface length is minimized to a value of $2\pi$.
			
\begin{figure}
		\begin{center}
			\hfill
			\subfigure[$t=0$]{
				\includegraphics[width=0.475\textwidth]{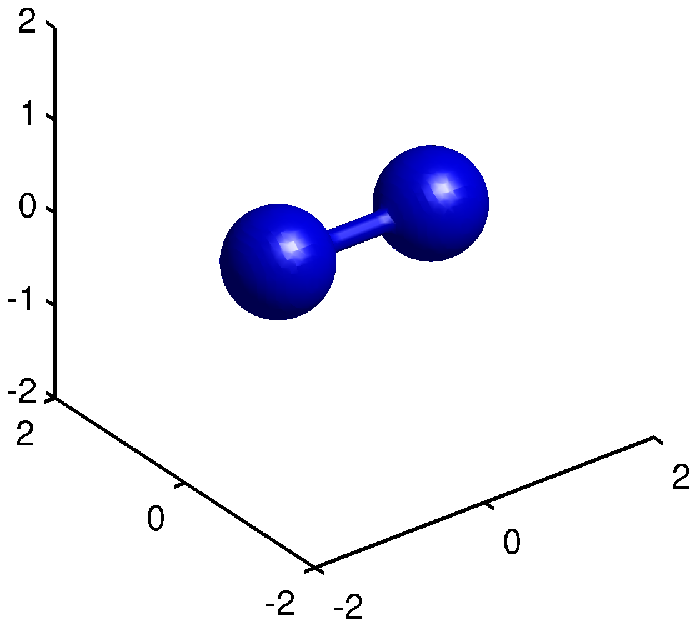}
				\label{fig:3d_dumbbell_a}
			}\hfill
			\subfigure[$t=0.0003$]{
				\includegraphics[width=0.475\textwidth]{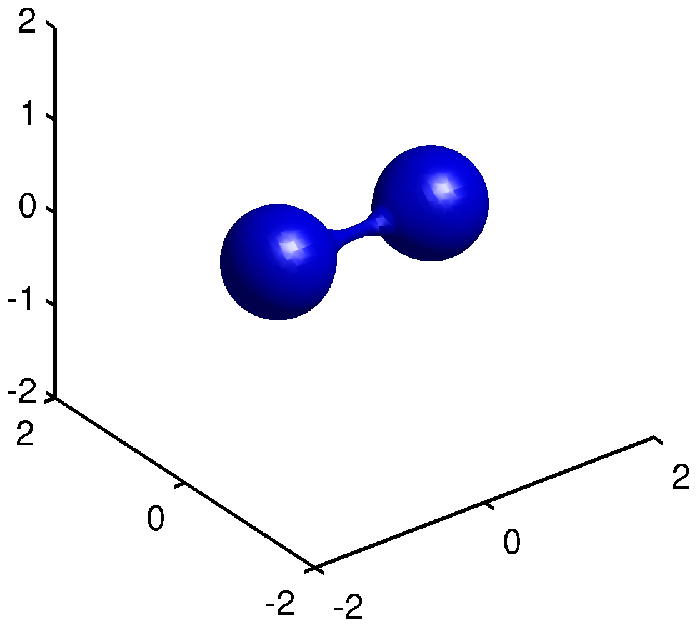}
				\label{fig:3d_dumbbell_b}
			} \hfill \\
			\hfill
			\subfigure [$t=0.0005$]{
				\includegraphics[width=0.475\textwidth]{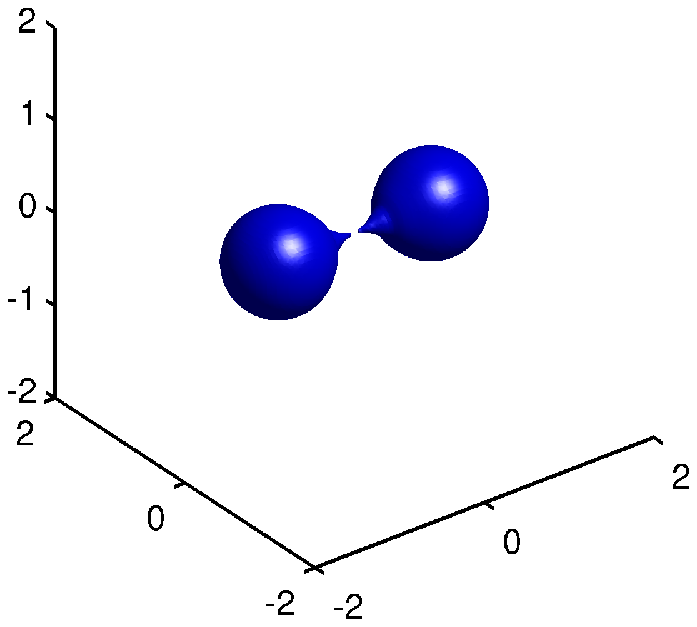}
				\label{fig:3d_dumbbell_c}
			}\hfill
			\subfigure [$t=0.005$]{
				\includegraphics[width=0.475\textwidth]{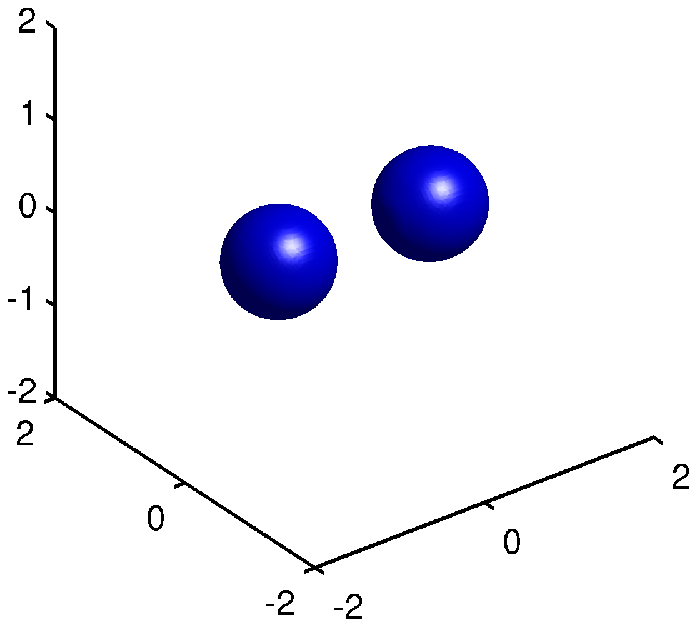}
				\label{fig:3d_dumbbell_d}
			}\hfill \\ 		
		\end{center}
		\caption{Motion by surface diffusion for a three-dimensional dumbbell surface. The grid spacing is $h=0.0625$ and the time step is $\Delta t=0.0001$}
		\label{fig:3d_dumbbell} 
	\end{figure}

\begin{figure}
		\begin{center}
			\hfill
			\subfigure[$t=0$]{
				\includegraphics[width=0.475\textwidth]{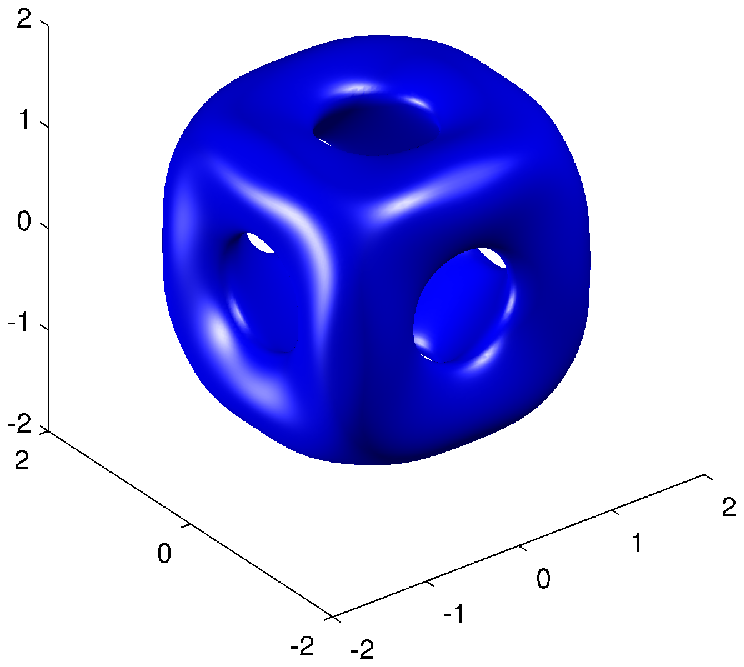}
				\label{fig:3d_complex_a}
			}\hfill
			\subfigure[$t=0.05$]{
				\includegraphics[width=0.475\textwidth]{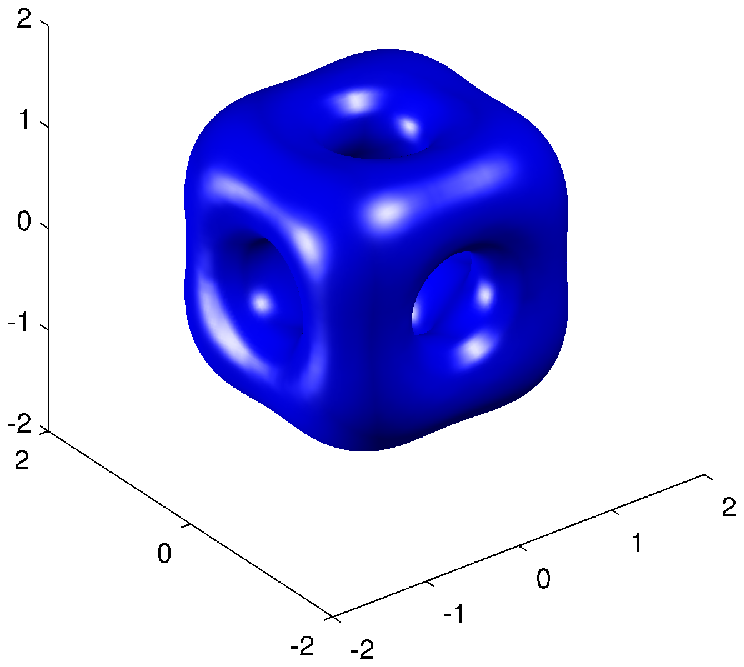}
				\label{fig:3d_complex_b}
			} \hfill \\
			\hfill
			\subfigure [$t=0.09$]{
				\includegraphics[width=0.475\textwidth]{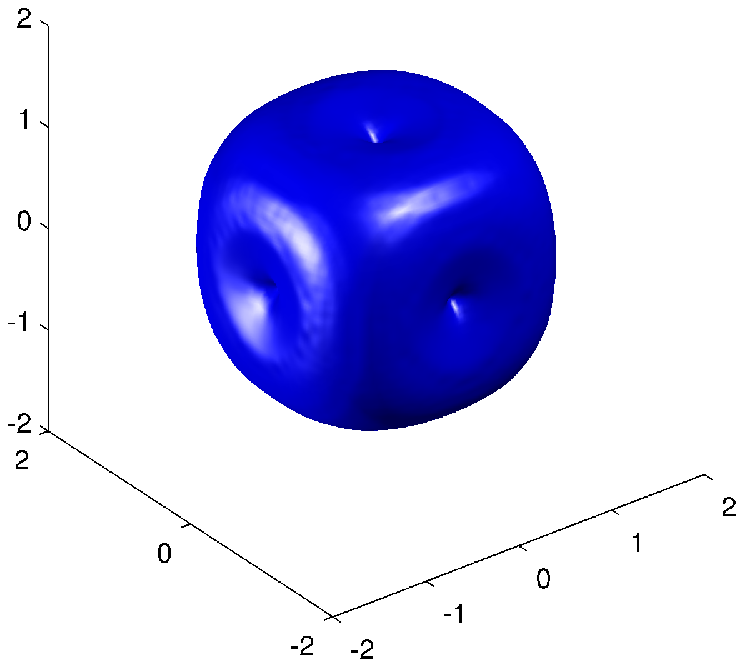}
				\label{fig:3d_complex_c}
			}\hfill
			\subfigure [$t=0.3$]{
				\includegraphics[width=0.475\textwidth]{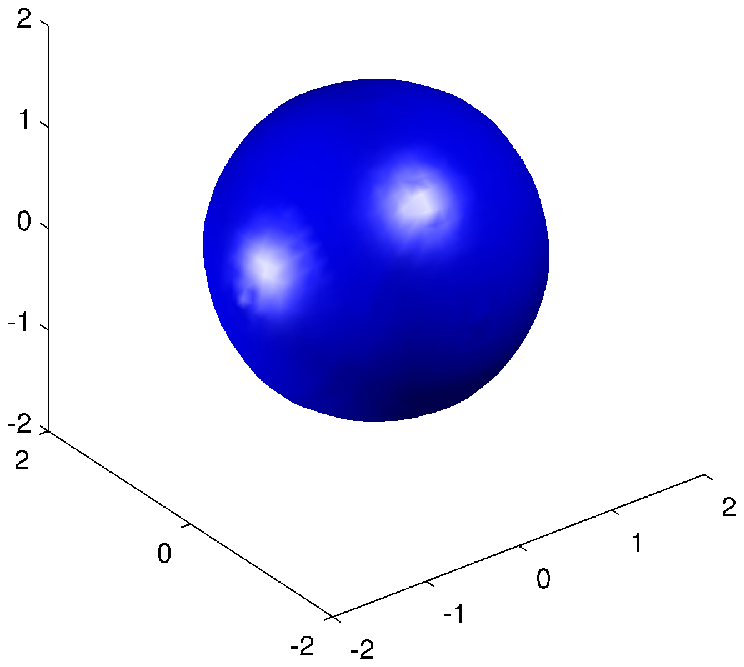}
				\label{fig:3d_complex_d}
			}\hfill \\ 		
		\end{center}
			   		
		\caption{Evolution of a complex three-dimensional shape under surface diffusion. The grid spacing is $h=0.0625$ and the time step is $\Delta t=10^{-4}$}
		\label{fig:3d_complex}												
	\end{figure}

	Next consider two three-dimensional surfaces, a dumbbell surface, Fig. \ref{fig:3d_dumbbell}, and a more complicated box-like surface, Fig. \ref{fig:3d_complex}.
	The dumbbell results in two distinct spheres due to the pinching-off of the thin region while the more complicated box-like surface results in a single sphere.

\begin{figure}
		\begin{center}

			\subfigure[$t=0$]{
  				\includegraphics[width=0.31\textwidth]{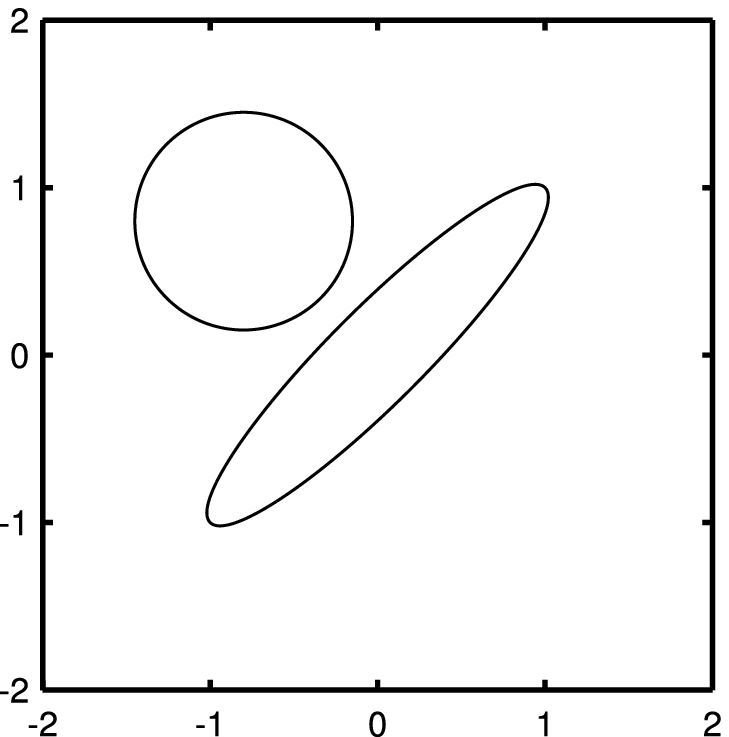}
			}
			\subfigure[$t=0.01$]{
  				\includegraphics[width=0.31\textwidth]{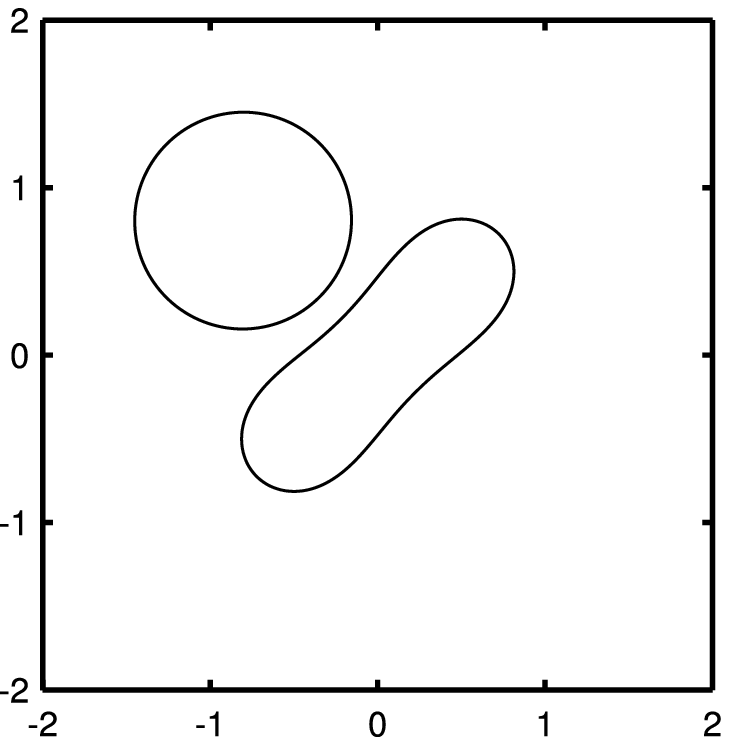}
			} 
			\subfigure[$t=0.029$] {
  				\includegraphics[width=0.31\textwidth]{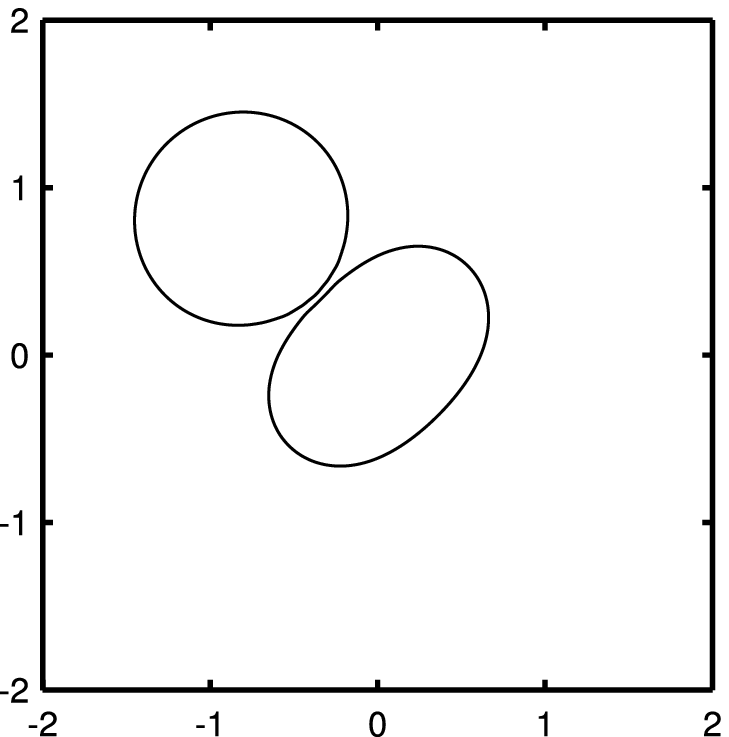}
			} \\
			\subfigure[$t=0.0304$] {
  				\includegraphics[width=0.31\textwidth]{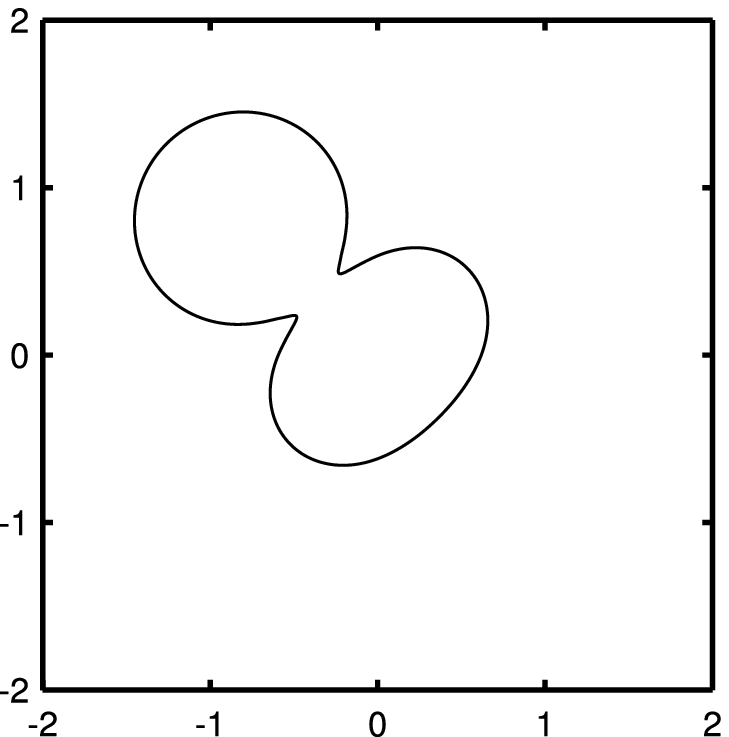}
			} 
			\subfigure[$t=0.032$] {
  				\includegraphics[width=0.31\textwidth]{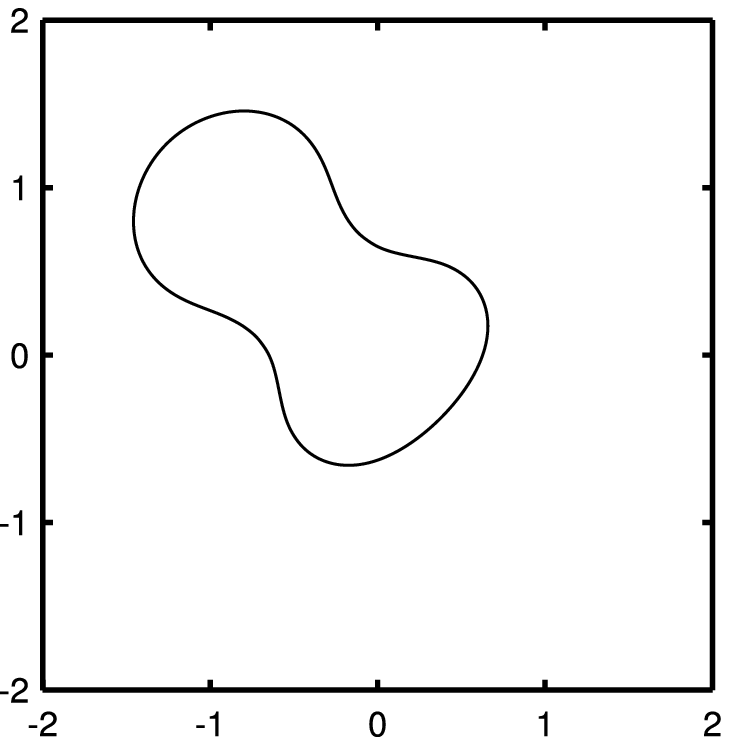}
			}
			\subfigure[$t=0.035$] {
  				\includegraphics[width=0.31\textwidth]{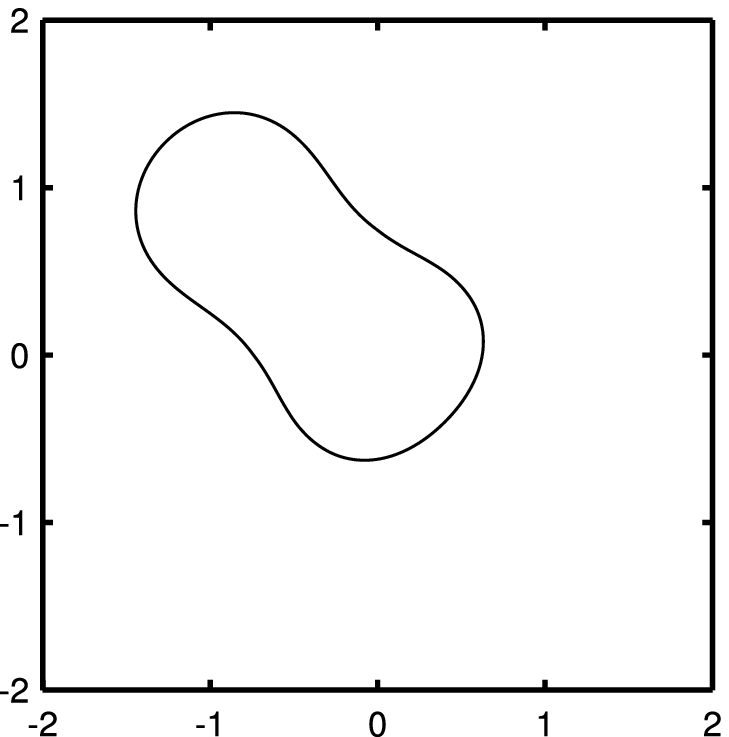}
			}\\  					
			\caption{Coalescence of two bodies under surface diffusion flow. The time step is $\Delta t=0.0001$ and while the grid spacing is $h=0.03125$. There is minimal
			spurious motion of the circle before merging.}	
			\label{fig:merging}
		\end{center}				
	\end{figure}

\begin{figure}
		\begin{center}
			\subfigure[$t=0.0235$] {
  				\includegraphics[width=0.475\textwidth]{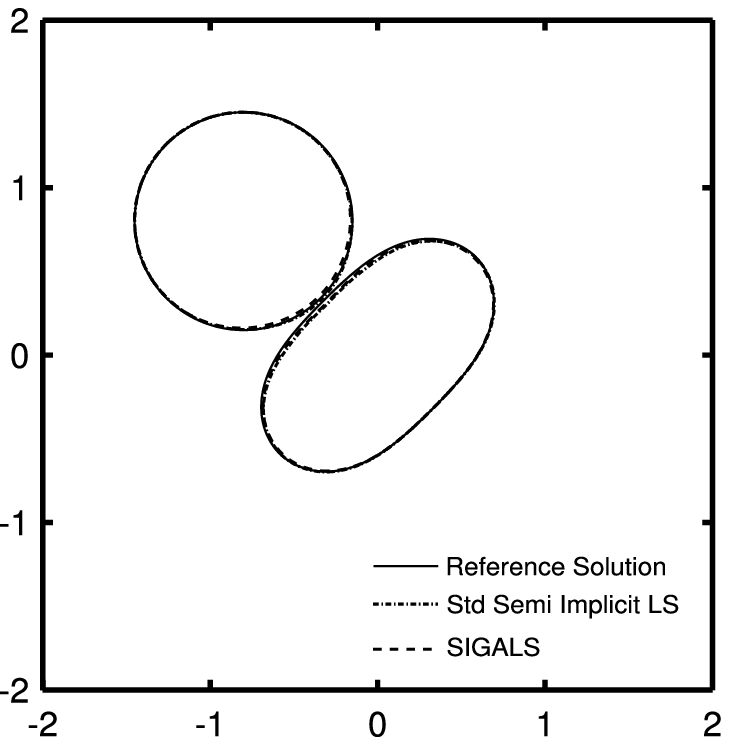}
				\label{fig:merging_comparison_a}
			} 			
			\subfigure[$t=0.025$]{
  				\includegraphics[width=0.475\textwidth]{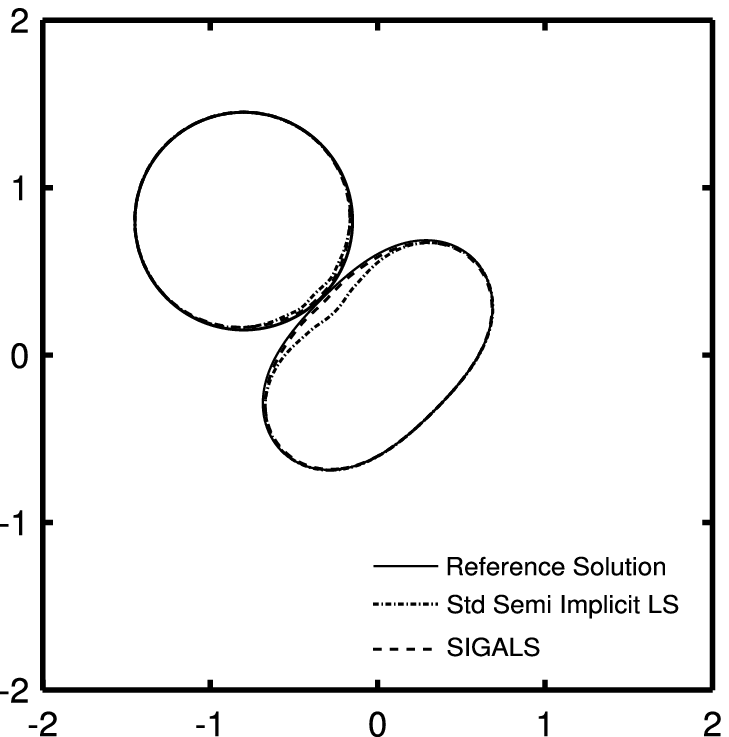}
				\label{fig:merging_comparison_b}
			}\\
			\subfigure[Circle] {
  				\includegraphics[width=0.475\textwidth]{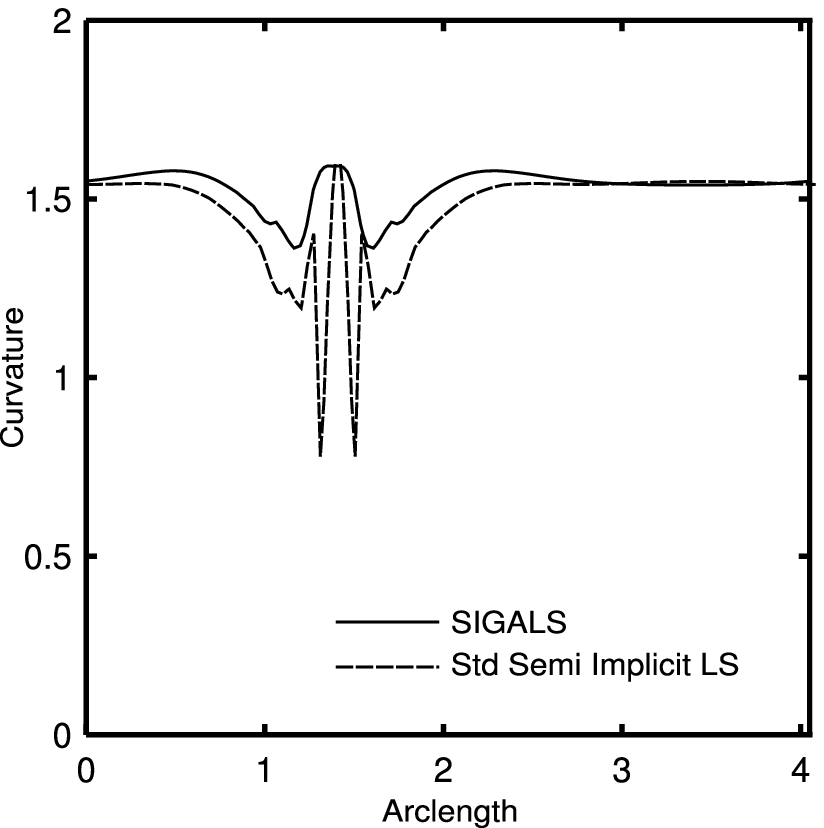}
				\label{fig:merging_comparison_c}
			}
			\subfigure[Ellipse] {
  				\includegraphics[width=0.475\textwidth]{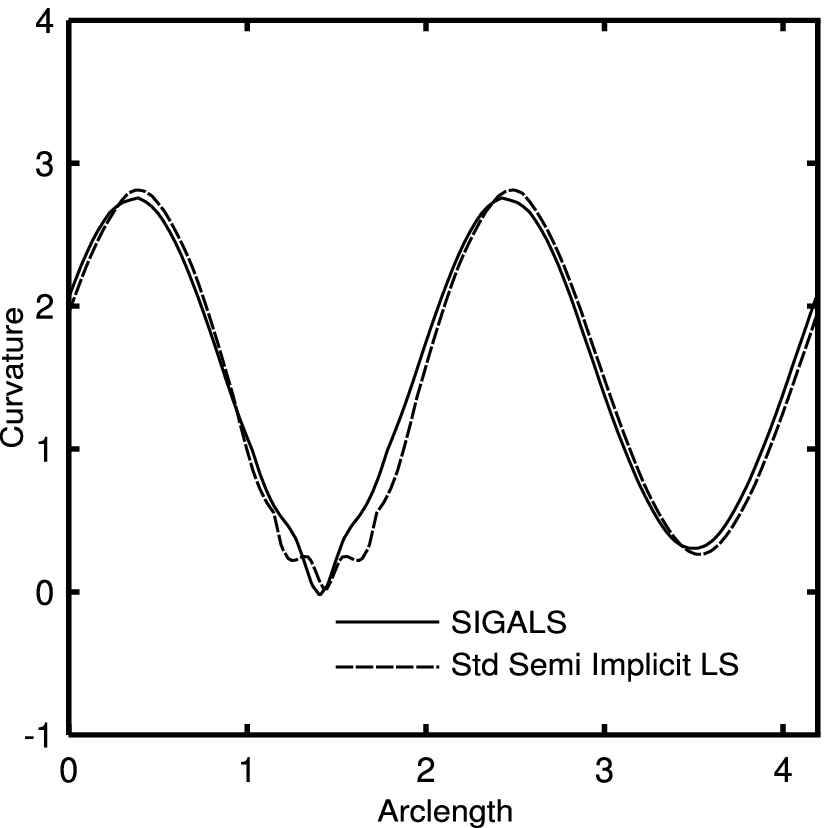}
				\label{fig:merging_comparison_d}
			}\\  					
			\caption{A comparison between the standard semi-implicit level set, the SIGALS method, and the reference solution for the case of a merging circle and ellipse.
			Two representative times are shown in (a) and (b). The curvature along the (c) circle and (d) ellipse are also shown for a time of $t=0.0235$. 
			The large curvature oscillations	in the standard semi-implicit level set result in the incorrect behavior observed in (b).
			}	
			\label{fig:merging_comparison}
		\end{center}
		
	\end{figure}

	Finally consider the merging of a circle and ellipse under surface diffusion, Fig. \ref{fig:merging}. 
	The circle should remain stationary until the evolving ellipse
	contacts it. 	
	In the original semi-implicit level set work it was demonstrated that the non-local nature of the smoothing operator results in slight perturbation of the
	interfaces \cite{Smereka2003}. In this work a similar result is noticed but to a smaller degree than in the standard semi-implicit scheme. This is further demonstrated in 
	Fig. \ref{fig:merging_comparison}, where a comparison between the standard semi-implicit level set solution, the SIGALS method, and the reference solution is shown
	for times of $t=0.0235$ and $t=0.025$.
	The reference solution was obtained by evolving the interfaces independently. While both the regular and gradient-augmented semi-implicit level set schemes appear to be
	working correctly at a time of $t=0.0235$, Fig. \ref{fig:merging_comparison_a}, the standard semi-implicit scheme clearly demonstrates incorrect behavior 
	at time $t=0.025$, Fig. \ref{fig:merging_comparison_b}. This is due to large oscillations in the calculated curvature, 
	shown in Figs. \ref{fig:merging_comparison_c} and \ref{fig:merging_comparison_d} for the circle and ellipse interfaces, respectively, 
	at the time $t=0.0235$. While oscillations still occur in the SIGALS method they are smaller than in standard method. This is due to the additional information 
	provided by explicitly tracking the gradient of the level set.

%%%%%%%%%%%%%%%%%%%%%%%%%%%%%%%%%%%%%%%%%%%%%%%%%%%%%%%%%%%%%%%%%%%%%%%%%%%%%%%%%%%%%%%%%%%%%%%%%%%%%%%%%%%%%%%%%%%%%%%%%%%	
\section{Summary}
\label{sec:4.0}
	This work represents a combination of the standard semi-implicit level set method and the gradient-augmented level set method. The result is a hybrid Lagrangian-Eulerian
	method that is easily applied in both two- and three-dimensions. 
	Sample results were provided for motion by mean curvature and surface diffusion. The addition of gradient information reduced the
	non-local nature of the semi-implicit scheme and allowed for investigations of flows depending on high-order derivatives of the interface.

\clearpage
\newpage

\bibliography{sigals}

\begin{thebibliography}{10}
\expandafter\ifx\csname url\endcsname\relax
  \def\url#1{\texttt{#1}}\fi
\expandafter\ifx\csname urlprefix\endcsname\relax\def\urlprefix{URL }\fi
\expandafter\ifx\csname href\endcsname\relax
  \def\href#1#2{#2} \def\path#1{#1}\fi

\bibitem{Li1999}
Z.~Li, H.~Zhao, H.~Gao, {A numerical study of electro-migration voiding by
  evolving level set functions on a fixed Cartesian grid}, Journal of
  Computational Physics {152} (1999) {281--304}.

\bibitem{Luo2008}
J.~Luo, Z.~Luo, L.~Chen, L.~Tong, M.~Wang, {A semi-implicit level set method
  for structural shape and topology optimization}, Journal of Computational
  Physics 227~(11) (2008) 5561--5581.

\bibitem{Malladi1995}
R.~Malladi, J.~Sethian, {Image processing via level set curvature flow.},
  Proceedings of the National Academy of Sciences of the United States of
  America 92~(15) (1995) 7046--50.

\bibitem{Sethian2003}
J.~A. Sethian, P.~Smereka, {Level Set Methods for Fluid Interfaces}, Annual
  Review of Fluid Mechanics 35~(1) (2003) 341--372.

\bibitem{Osher1988}
S.~Osher, J.~Sethian, Fronts propagating with curvature-dependent speed -
  algorithms based on hamilton-jacobi formulations, Journal of Computational
  Physics 79~(1) (1988) 12--49.

\bibitem{Sussman1994}
M.~Sussman, {A Level Set Approach for Computing Solutions to Incompressible
  Two-Phase Flow} (Sep. 1994).

\bibitem{Salac2011}
D.~Salac, M.~Miksis, A level set projection model of lipid vesicles in general
  flows, Journal of Computational Physics 230~(22) (2011) 8192 -- 8215.

\bibitem{Adalsteinsson}
D.~Adalsteinsson, J.~Sethian, {A level set approach to a unified model for
  etching, deposition, and lithography .3. Redeposition, reemission, surface
  diffusion, and complex simulations}, Journal of Computational Physics {138}
  (197) {193--223}.

\bibitem{MERRIMAN1994}
B.~Merriman, J.~Bence, S.~Osher, Motion of multple junctions - a level set
  approach, Journal of Computational Physics {112} (1994) {334--363}.

\bibitem{Mullins1957}
W.~W. Mullins, Theory of thermal grooving, Journal of Applied Physics 28~(3)
  (1957) 333--339.

\bibitem{Veerapaneni2009}
S.~Veerapaneni, D.~Gueyffier, D.~Zorin, G.~Biros, A boundary integral method
  for simulating the dynamics of inextensible vesicles suspended in a viscous
  fluid in 2d, Journal of Computational Physics 228~(7) (2009) 2334--2353.

\bibitem{chopp1999}
D.~Chopp, J.~Sethian, {Motion by Intrinsic Laplacian of Curvature}, Interfaces
  and Free Boundaries 1 (1999) 1--18.

\bibitem{Khenner2001}
M.~Khenner, A.~Averbuch, M.~Israeli, M.~Nathan, {Numerical Simulation of Grain
  Boundary Grooving By Level Set Method }, Journal of Computational Physics
  (2001) 764--784.

\bibitem{Salac2008}
D.~Salac, W.~Lu, {A Local Semi-Implicit Level-Set Method for Interface Motion},
  Journal of Scientific Computing 35~(2-3) (2008) 330--349.

\bibitem{Smereka2003}
P.~Smereka, Semi-implicit level set methods for curvature and surface diffusion
  motion, Journal of Scientific Computing 19 (2003) 439--456.

\bibitem{Nave2010}
J.-C. Nave, R.~R. Rosales, B.~Seibold, {A gradient-augmented level set method
  with an optimally local, coherent advection scheme} (May 2010).

\bibitem{Peng1999}
D.~Peng, {A PDE-Based Fast Local Level Set Method}, Journal of Computational
  Physics 155~(2) (1999) 410--438.

\bibitem{Kumar2004}
A.~Kumar, {Isotropic finite-differences}, Journal of Computational Physics
  201~(1) (2004) 109--118.

\bibitem{Vorst1992}
H.~A. van~der Vorst, Bi-cgstab: A fast and smoothly converging variant of bi-cg
  for the solution of nonsymmetric linear systems, SIAM Journal on Scientific
  Computing 13~(2) (1992) 631--644.

\bibitem{Saad2003}
Y.~Saad, {Iterative Methods for Sparse Linear Systems, Second Edition}, 2nd
  Edition, Society for Industrial and Applied Mathematics, 2003.

\end{thebibliography}

\end{document}